\documentclass{elsart}

\usepackage[T1]{fontenc}
\usepackage[latin1]{inputenc}
\usepackage{latexsym,color}
\usepackage{amsmath}
\usepackage{natbib}
\usepackage{stmaryrd}
\usepackage{subfigure}
\usepackage{amssymb}
\usepackage{amsfonts}
\usepackage{times,verbatim}
\usepackage{graphicx,lscape}
\usepackage{psfrag}
\usepackage{subfigure}
\usepackage{epsfig}
\usepackage{array}
\usepackage{flafter}
\usepackage{bm}
\usepackage[named]{algo}
\usepackage[all,2cell,ps]{xy} \UseAllTwocells \SilentMatrices
\usepackage{fancyheadings}
\usepackage{multicol}

\newcommand{\expect}{\mathop{\mathbb{E}}}
\newcommand{\Prob}{\mathrm{P}}
\newcommand{\GO}{\mathrm{O}}
\newcommand{\Tr}{^{\sf T}}
\newcommand{\argmin}{\mathop{\mathrm{arg\,min}}}
\newcommand{\minim}{\mathop{\mathrm{min}}}
\newcommand{\var}{\mathop{\mathrm{var}}}
\newcommand{\grid}{\mathbb{G}}
\newcommand{\DD}{\mathbb{D}}
\newcommand{\Sample}{\mathbb{S}}
\newcommand{\XX}{\mathbb{X}}
\newcommand{\HH}{\mathbb{H}}
\newcommand{\minimizer}{\bm{x}^{\ast}}
\newcommand{\Minimizer}{\bm{X}^{\ast}}
\newcommand{\Minimizergrid}{\Minimizer_{\grid}}
\newcommand{\pred}{\hat{F}}
\newcommand{\x}{\bm{x}}
\newcommand{\y}{\bm{y}}
\newcommand{\lambdas}{\bm{\lambda}}
\newcommand{\covmat}{\bm{K}}
\newcommand{\predvar}{\hat{\sigma}}
\newcommand{\predmean}{\hat{f}}
\newcommand{\param}{\bm{\theta}}
\newcommand{\ngrid}{N}
\newcommand{\nobs}{n}
\newcommand{\ntraj}{r}
\newcommand{\npart}{M}

\newcommand{\ddpcond}{p_{\Minimizer_{\grid}|\bm{f_\Sample}}} 
\newcommand{\ddpaprox}{\hat{p}_{\Minimizer_{\grid}|\bm{f_\Sample}}} 

\newcommand{\quant}{Q^\alpha}
\newcommand{\event}{\mathcal{B}}
\newcommand{\MinimizerSet}{\mathcal{M}^*}

\begin{document}

\maketitle

\begin{frontmatter}

\title{An informational approach to the global optimization
  of expensive-to-evaluate functions}

\author{Julien Villemonteix$^{*\dagger\ddagger}$, 
  Emmanuel Vazquez$^{\dagger}$ and Eric Walter$^{*}$}
 
\address{$^\dagger$ Sup\'elec, D\'epartement Signaux et Syst\`emes \'Electroniques  \\
         91192 Gif-sur-Yvette, France \\
         $^*$ Laboratoire des Signaux et Syst\`emes,
         CNRS-Sup\'elec-Univ Paris-Sud \\
         $\ddagger$ Renault S.A., service 68240 \\
         78298 Guyancourt, France}

\begin{abstract}
  In many global optimization problems motivated by engineering
  applications, the number of function evaluations is severely limited
  by time or cost. To ensure that each evaluation
  contributes to the localization of good candidates for the role of
  global minimizer, a sequential choice of evaluation points is
  usually carried out. In particular, when Kriging is used to
  interpolate past evaluations, the uncertainty associated with the
  lack of information on the function can be expressed and used to
  compute a number of criteria accounting for the interest of an
  additional evaluation at any given point. This paper introduces
  minimizer entropy as a new Kriging-based criterion for the
  sequential choice of points at which the function should be
  evaluated. Based on \emph{stepwise uncertainty reduction}, it
  accounts for the informational gain on the minimizer expected from a
  new evaluation.  The criterion is approximated using conditional
  simulations of the Gaussian process model behind Kriging, and then
  inserted into an algorithm similar in spirit to the \emph{Efficient Global
    Optimization} (EGO) algorithm. An empirical comparison is carried
  out between our criterion and \emph{expected improvement}, one of
  the reference criteria in the literature. Experimental results indicate
  major evaluation savings over EGO. Finally, the method, which we call
  IAGO (for Informational Approach to Global Optimization) is extended
  to robust optimization problems, where both the factors to be tuned
  and the function evaluations are corrupted by noise.
\end{abstract}


\end{frontmatter}

\section{Introduction}

This paper is devoted to global optimization in a context of expensive
function evaluation. The objective is to find global minimizers in
$\XX$ (the factor space, a bounded subset of $\mathbb{R}^d$) of an
unknown function $f :\XX \to \mathbb{R} $, using a very limited number
of function evaluations. Note that the global minimizer may not be
unique (any global minimizer will be denoted as $\minimizer$).
Such a problem is frequently encountered in the industrial world. For
instance, in the automotive industry, optimal crash-related parameters
are obtained using costly real tests and time-consuming computer
simulations (a single simulation of crash-related deformations may
take up to 24 hours on dedicated servers). It then becomes essential
to favor optimization methods that use the dramatically scarce
information as efficiently as possible.

To make up for the lack of knowledge on the function, surrogate (also
called meta or approximate) models are used to obtain cheap
approximations (\cite{Jones2001}). They turn out to be convenient
tools for visualizing the function behavior or suggesting the location
of an additional point at which $f$ should be evaluated in the search
for $\minimizer$.  Surrogate models based on Gaussian processes have
received particular attention. Known in Geostatistics under the name
of \emph{Kriging} since the early 1960s (\cite{Matheron1963}),
Gaussian process models provide a probabilistic framework to account
for the uncertainty stemming from the lack of information on the
system. When dealing with an optimization problem, this framework
allows the set of function evaluations to be chosen efficiently
(\cite{Jones2001, Jones1998, Huang2006}).

In this context, several strategies have been proposed, with
significant advantages over traditional optimization methods when
confronted to expensive-to-evaluate functions. Most of them
\emph{implicitly} seek a likely value for $\minimizer$, and
then assume it to be a suitable location for a new evaluation of $f$.
Yet, given existing evaluation results, the most likely location of a global
minimizer is not necessarily a good evaluation point to improve our
knowledge on $\minimizer$. As we shall show, by making
full use of Kriging, it is instead possible to \emph{explicitly}
estimate the probability distribution of the optimum location, which
allows an information-based search strategy.

Based on these observations, the present paper introduces minimizer
entropy as a criterion for the choice of new evaluation points. This
criterion, directly inspired from \emph{stepwise uncertainty reduction}
(\cite{Geaman1995}), is then inserted in an algorithm similar to the
\emph{Efficient Global Optimization} (EGO) algorithm
(\cite{Jones1998}). We call the resulting algorithm IAGO, for \emph{Informational
Approach to Global Optimization}.

Section~\ref{sec:adapt-design-crit} recalls the principle of
Kriging-based optimization, along with some general ideas on Gaussian
process modeling that are used in Section~\ref{sec:krig-cond-simul} to
build an estimate of the distribution of the global minimizers.
Section~\ref{sec:stepw-uncert-reduct} details the stepwise uncertainty
reduction approach applied to global optimization, while
Section~\ref{sec:impl-sur-strat} describes the corresponding algorithm
and its extensions to noisy problems. Section~\ref{sec:illustrations}
illustrates the behavior of the new algorithm on some simple benchmark
problems, along with its performances compared with those of the
classical EGO algorithm, chosen for its good compromise between local
and global search (\cite{Sasena2002}).  Finally, after a conclusion
section and to make this paper self-contained,
Section~\ref{sec:annex-model-with} presents, as an appendix, some more
results on Gaussian process modeling and Kriging.
\section{Kriging-based global optimization}
\label{sec:adapt-design-crit}
When dealing with expensive-to-evaluate functions, optimization
methods based on probabilistic surrogate models (and Kriging in
particular) have significant advantages over traditional optimization
techniques, as they require fewer function evaluations. Kriging can
indeed provide a cheap and accurate approximation of the function, but
also an estimate of the potential error in this approximation.
Numerous illustrations of this superiority can be found in the
literature (see, for instance, \cite{cox1997}) and many variations
have been explored (for extensive surveys, see \cite{Jones2001} and
\cite{Sasena2002}). As explained in this section, these methods deal
with the cost of evaluation using an adaptive sampling strategy,
replacing the optimization of the expensive-to-evaluate function $f$ by a
series of optimizations of a cheap criterion.
\subsection{Gaussian process modeling and Kriging}
\label{sec:gauss-proc-model}
This section briefly recalls the principle of Gaussian process (GP)
modeling, and lays down the necessary notation. A more detailed
presentation is available in the appendix
(Section~\ref{sec:annex-model-with}).

When modeling with Gaussian processes, the function $f$ is assumed to
be a sample path of a second-order Gaussian random process $F$. If we
denote $(\Omega,\mathcal{A},\mathcal{P})$ the underlying probability space,
this amounts to assuming that $\exists\, \omega\in\Omega$, such that
$F(\omega,\cdot)=f(\cdot)$. Whenever possible, we shall omit the dependence of $F$
in $\omega$ to simplify notation. 

In particular, given a set of $\nobs$ evaluation points
$\Sample=\{\x_1,\ldots,\x_{\nobs}\}$ (the \emph{design}), $\forall \x_{i}\in \Sample$
the evaluation result $f(\x_i)$ is viewed as a sample path of the
random variable $F(\x_i)$.  Kriging computes an unbiased linear
predictor of $F(\x)$ in the vector space
$\HH_\Sample=\mathrm{span}\{F(\x_1),\ldots,F(\x_n)\}$, which can be written
as
\begin{equation}
\label{eq:prediction}
   \pred(\x) = \lambdas(\x)\Tr\bm{F}_\Sample\,,
\end{equation}
with $\bm{F}_\Sample=[F(\x_1),\ldots,F(\x_n)]\Tr$, and $\lambdas(\x)$ the
vector of Kriging coefficients for the prediction at $\x$.

Given the covariance of $F$, the Kriging coefficients can be computed
along with the variance of the prediction error
\begin{equation}
\predvar^2(\x)=\var(\pred(\x)-F(\x)).
\label{eq:MSEsimple}
\end{equation}
The covariance of $F$ is chosen within a class of parametrized covariances
(for instance, the Mat\`ern class), and its parameters are either
estimated from the data or chosen a priori (see
Section~\ref{sec:estim-param} for details on the choice of a
covariance).

Once $f$ has been evaluated at all evaluation points in $\Sample$, the
predicted value of $f$ at $\x$ is given by
\begin{equation}
  \predmean (\x)=\lambdas(\x)\Tr\bm{f}_\Sample\:,
\label{eq:Pred}
\end{equation}
with $\bm{f}_{\Sample}=[f(\x_1),\ldots,f(\x_n)]\Tr$ ($\bm{f}_{\Sample}$ is
viewed as a sample value of $\bm{F}_{\Sample}$). The same results could
be derived in a Bayesian framework, where $F(\x)$ is Gaussian
conditionally to the evaluations carried out
($\bm{F}_{\Sample}=\bm{f}_{\Sample}$), with mean $\predmean(\x)$ and
variance $\predvar^2(\x)$.

Note that
\begin{equation}
  \label{eq:prediction_at_evaluation_points}
    \forall \:\x_i \:\in\: \Sample,\:\: \pred(\x_i)=F(\x_i), 
\end{equation}
and that the prediction at $\x_i \:\in\: \Sample$ is $f(\x_i)$. When $f$
is assumed to be evaluated exactly, Kriging is thus an interpolation,
with the considerable advantage over other interpolation methods that it
also provides an explicit characterization of the prediction error
(zero-mean Gaussian with variance $\predvar^2(\x)$).
\subsection{Adaptive sampling strategies}
\begin{figure}
  \centering
  % This file is generated by the MATLAB m-file laprint.m. It can be included
% into LaTeX documents using the packages graphicx, color and psfrag.
% It is accompanied by a postscript file. A sample LaTeX file is:
%    \documentclass{article}\usepackage{graphicx,color,psfrag}
%    \begin{document}\input{optimsimple}\end{document}
% See http://www.mathworks.de/matlabcentral/fileexchange/loadFile.do?objectId=4638
% for recent versions of laprint.m.
%
% created by:           LaPrint version 3.16 (13.9.2004)
% created on:           20-Jun-2006 16:40:22
% eps bounding box:     15 cm x 10.7813 cm
% comment:              
%
\begin{psfrags}%
\psfragscanon%
%
% text strings:
\psfrag{s03}[b][b]{\color[rgb]{0,0,0}\setlength{\tabcolsep}{0pt}\begin{tabular}{c}$f$,$\predmean$\end{tabular}}%
\psfrag{s05}[t][t]{\color[rgb]{0,0,0}\setlength{\tabcolsep}{0pt}\begin{tabular}{c}$\x$\end{tabular}}%
\psfrag{s06}[b][b]{\color[rgb]{0,0,0}\setlength{\tabcolsep}{0pt}\begin{tabular}{c}$f$,$\predmean$\end{tabular}}%
%
% xticklabels:
\psfrag{x01}[t][t]{0}%
\psfrag{x02}[t][t]{0.1}%
\psfrag{x03}[t][t]{0.2}%
\psfrag{x04}[t][t]{0.3}%
\psfrag{x05}[t][t]{0.4}%
\psfrag{x06}[t][t]{0.5}%
\psfrag{x07}[t][t]{0.6}%
\psfrag{x08}[t][t]{0.7}%
\psfrag{x09}[t][t]{0.8}%
\psfrag{x10}[t][t]{0.9}%
\psfrag{x11}[t][t]{1}%
\psfrag{x12}[t][t]{0}%
\psfrag{x13}[t][t]{0.1}%
\psfrag{x14}[t][t]{0.2}%
\psfrag{x15}[t][t]{0.3}%
\psfrag{x16}[t][t]{0.4}%
\psfrag{x17}[t][t]{0.5}%
\psfrag{x18}[t][t]{0.6}%
\psfrag{x19}[t][t]{0.7}%
\psfrag{x20}[t][t]{0.8}%
\psfrag{x21}[t][t]{0.9}%
\psfrag{x22}[t][t]{1}%
\psfrag{x23}[t][t]{2}%
\psfrag{x24}[t][t]{4}%
\psfrag{x25}[t][t]{6}%
\psfrag{x26}[t][t]{8}%
\psfrag{x27}[t][t]{10}%
\psfrag{x28}[t][t]{12}%
\psfrag{x29}[t][t]{2}%
\psfrag{x30}[t][t]{4}%
\psfrag{x31}[t][t]{6}%
\psfrag{x32}[t][t]{8}%
\psfrag{x33}[t][t]{10}%
\psfrag{x34}[t][t]{12}%
%
% yticklabels:
\psfrag{v01}[r][r]{0}%
\psfrag{v02}[r][r]{0.1}%
\psfrag{v03}[r][r]{0.2}%
\psfrag{v04}[r][r]{0.3}%
\psfrag{v05}[r][r]{0.4}%
\psfrag{v06}[r][r]{0.5}%
\psfrag{v07}[r][r]{0.6}%
\psfrag{v08}[r][r]{0.7}%
\psfrag{v09}[r][r]{0.8}%
\psfrag{v10}[r][r]{0.9}%
\psfrag{v11}[r][r]{1}%
\psfrag{v12}[r][r]{0}%
\psfrag{v13}[r][r]{0.1}%
\psfrag{v14}[r][r]{0.2}%
\psfrag{v15}[r][r]{0.3}%
\psfrag{v16}[r][r]{0.4}%
\psfrag{v17}[r][r]{0.5}%
\psfrag{v18}[r][r]{0.6}%
\psfrag{v19}[r][r]{0.7}%
\psfrag{v20}[r][r]{0.8}%
\psfrag{v21}[r][r]{0.9}%
\psfrag{v22}[r][r]{1}%
\psfrag{v23}[r][r]{0}%
\psfrag{v24}[r][r]{2}%
\psfrag{v25}[r][r]{4}%
\psfrag{v26}[r][r]{6}%
\psfrag{v27}[r][r]{0}%
\psfrag{v28}[r][r]{2}%
\psfrag{v29}[r][r]{4}%
\psfrag{v30}[r][r]{6}%
%
% Figure:
\resizebox{12cm}{!}{\includegraphics{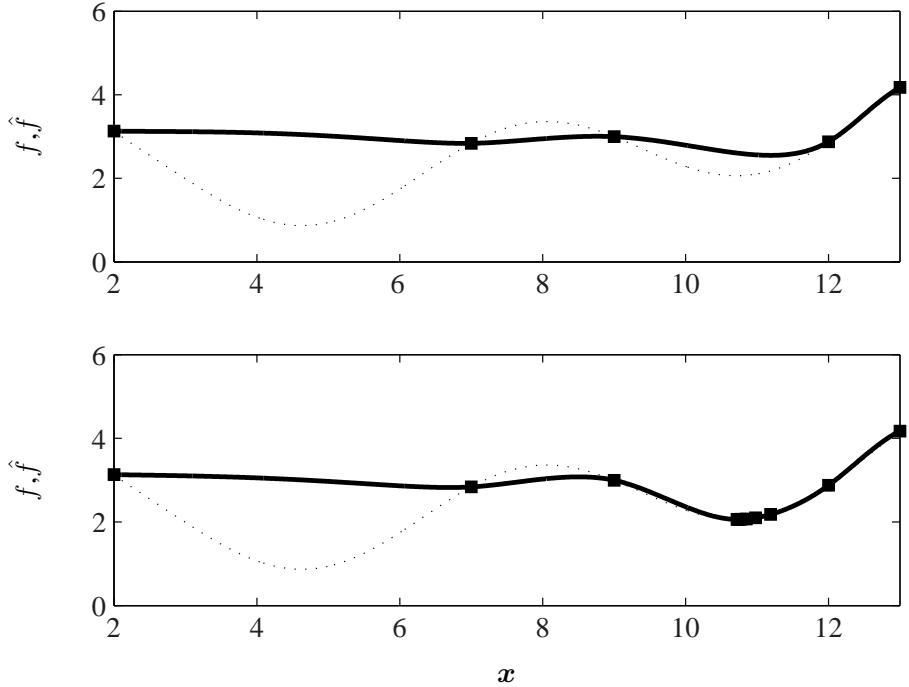}}%
\end{psfrags}%
%
% End optimsimple.tex

  \caption{Naive approach to optimization using Kriging: (\emph{top}) prediction
    $\predmean$ (bold line) of the true function $f$ (dotted line, supposedly
    unknown) obtained from an initial design materialized by squares;
    (\emph{bottom}) prediction after seven iterations of a direct
    minimization of $\predmean$.}
  \label{fig:optimsimple}
\end{figure}

The general principle of optimization using Kriging is iteratively to
evaluate $f$ at a point that optimizes a criterion based on the model
obtained using previous evaluation results. The simplest approach would
be to choose a minimizer of the prediction $\predmean$ as a new
evaluation point. However, by doing so, too much confidence would be put in
the current prediction and search is likely to stall on a local
optimum (as illustrated by Figure~\ref{fig:optimsimple}). To
compromise between local and global search, more emphasis has to be
put on the prediction error that can indicate locations where
additional evaluations are needed to improve confidence in the
model. This approach has led to a number of criteria, based on both
prediction and prediction error, designed to select additional
evaluation points.

\begin{figure}
  \centering
  % This file is generated by the MATLAB m-file laprint.m. It can be included
% into LaTeX documents using the packages graphicx, color and psfrag.
% It is accompanied by a postscript file. A sample LaTeX file is:
%    \documentclass{article}\usepackage{graphicx,color,psfrag}
%    \begin{document}\input{ei}\end{document}
% See http://www.mathworks.de/matlabcentral/fileexchange/loadFile.do?objectId=4638
% for recent versions of laprint.m.
%
% created by:           LaPrint version 3.16 (13.9.2004)
% created on:           01-Jun-2006 13:51:44
% eps bounding box:     15 cm x 10.7813 cm
% comment:              
%
\begin{psfrags}%
\psfragscanon%
%
% text strings:
\psfrag{s01}[t][t]{\color[rgb]{0,0,0}\setlength{\tabcolsep}{0pt}\begin{tabular}{c}\end{tabular}}%
\psfrag{s02}[b][b]{\color[rgb]{0,0,0}\setlength{\tabcolsep}{0pt}\begin{tabular}{c}$\predmean(\x)$\end{tabular}}%
\psfrag{s05}[t][t]{\color[rgb]{0,0,0}\setlength{\tabcolsep}{0pt}\begin{tabular}{c}$\x$\end{tabular}}%
\psfrag{s06}[b][b]{\color[rgb]{0,0,0}\setlength{\tabcolsep}{0pt}\begin{tabular}{c}$\mathrm{EI}(\x)$\end{tabular}}%
\psfrag{s17}[l][l]{\color[rgb]{0,0,0}\setlength{\tabcolsep}{0pt}\begin{tabular}{l}$\x_{\mathrm{new}}$ \end{tabular}}%
%
% xticklabels:
\psfrag{x01}[t][t]{0}%
\psfrag{x02}[t][t]{0.1}%
\psfrag{x03}[t][t]{0.2}%
\psfrag{x04}[t][t]{0.3}%
\psfrag{x05}[t][t]{0.4}%
\psfrag{x06}[t][t]{0.5}%
\psfrag{x07}[t][t]{0.6}%
\psfrag{x08}[t][t]{0.7}%
\psfrag{x09}[t][t]{0.8}%
\psfrag{x10}[t][t]{0.9}%
\psfrag{x11}[t][t]{1}%
\psfrag{x12}[t][t]{0}%
\psfrag{x13}[t][t]{0.1}%
\psfrag{x14}[t][t]{0.2}%
\psfrag{x15}[t][t]{0.3}%
\psfrag{x16}[t][t]{0.4}%
\psfrag{x17}[t][t]{0.5}%
\psfrag{x18}[t][t]{0.6}%
\psfrag{x19}[t][t]{0.7}%
\psfrag{x20}[t][t]{0.8}%
\psfrag{x21}[t][t]{0.9}%
\psfrag{x22}[t][t]{1}%
\psfrag{x23}[t][t]{2}%
\psfrag{x24}[t][t]{4}%
\psfrag{x25}[t][t]{6}%
\psfrag{x26}[t][t]{8}%
\psfrag{x27}[t][t]{10}%
\psfrag{x28}[t][t]{12}%
\psfrag{x29}[t][t]{2}%
\psfrag{x30}[t][t]{4}%
\psfrag{x31}[t][t]{6}%
\psfrag{x32}[t][t]{8}%
\psfrag{x33}[t][t]{10}%
\psfrag{x34}[t][t]{12}%
%
% yticklabels:
\psfrag{v01}[r][r]{0}%
\psfrag{v02}[r][r]{0.1}%
\psfrag{v03}[r][r]{0.2}%
\psfrag{v04}[r][r]{0.3}%
\psfrag{v05}[r][r]{0.4}%
\psfrag{v06}[r][r]{0.5}%
\psfrag{v07}[r][r]{0.6}%
\psfrag{v08}[r][r]{0.7}%
\psfrag{v09}[r][r]{0.8}%
\psfrag{v10}[r][r]{0.9}%
\psfrag{v11}[r][r]{1}%
\psfrag{v12}[r][r]{0}%
\psfrag{v13}[r][r]{0.1}%
\psfrag{v14}[r][r]{0.2}%
\psfrag{v15}[r][r]{0.3}%
\psfrag{v16}[r][r]{0.4}%
\psfrag{v17}[r][r]{0.5}%
\psfrag{v18}[r][r]{0.6}%
\psfrag{v19}[r][r]{0.7}%
\psfrag{v20}[r][r]{0.8}%
\psfrag{v21}[r][r]{0.9}%
\psfrag{v22}[r][r]{1}%
\psfrag{v23}[r][r]{0}%
\psfrag{v24}[r][r]{0.04}%
\psfrag{v25}[r][r]{0.08}%
\psfrag{v26}[r][r]{0}%
\psfrag{v27}[r][r]{2}%
\psfrag{v28}[r][r]{4}%
\psfrag{v29}[r][r]{6}%
%
% Figure:
\resizebox{12cm}{!}{\includegraphics{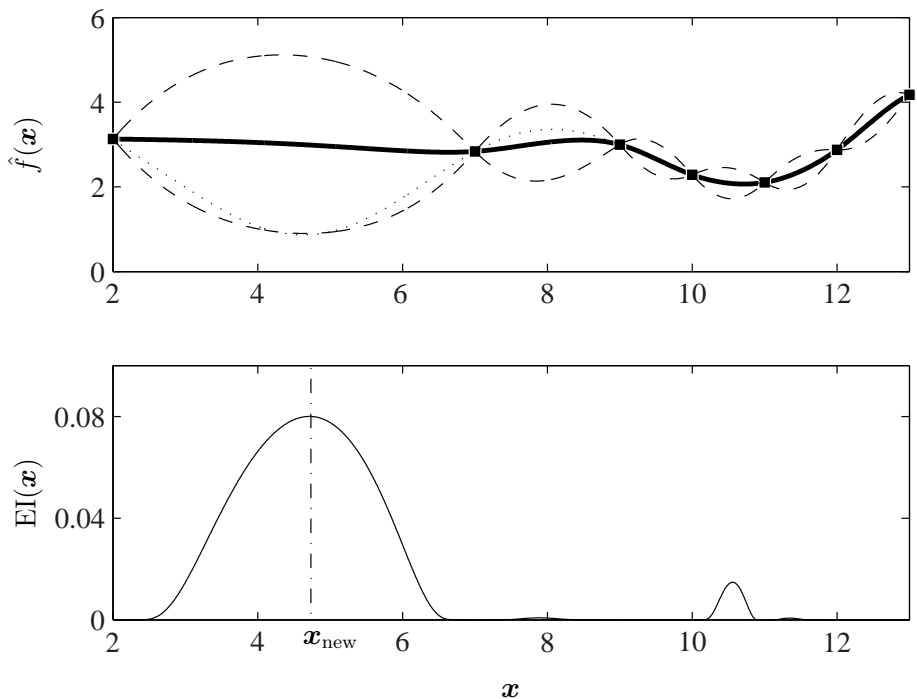}}%
\end{psfrags}%
%
% End ei.tex
 
  \caption{EI approach to optimization using
    Kriging: (\emph{top}) $\predmean$ (bold line), 95\% confidence intervals
    computed using $\predvar$ (dashed line) and true function $f$
    (dotted line); (\emph{bottom}) expected improvement. }
  \label{fig:EI}
\end{figure}

A standard example of such a criterion is \emph{expected
  improvement} (EI). As the name suggests, it involves computing how
much improvement in the optimum is expected, if $f$ is evaluated at a
given additional point. Let $F$ be the Gaussian process model, as
before, and $f_{\mathrm{min}}$ be the current best function value
obtained.  The improvement expected from an additional evaluation of
$f$ at $\x$ given $\bm{f}_\Sample$, the results of past evaluations, can
then be expressed as
$$
\mathrm{EI}(\x)=\expect\left[\max\left(f_{\mathrm{min}}-F\left(\x\right),0\right)|\bm{F}_\Sample\,=\,\bm{f}_\Sample\right].
$$    
Since $F(\x)$ is conditionally Gaussian with mean $\predmean(\x)$ and
variance $\predvar^2(\x)$, a convenient expression appears
(\cite{Schonlau1997}):
\begin{equation}
  \label{eq:EI}
  \mathrm{EI}(\x)=\predvar(\x)\left[ u\Phi(u)+ \frac{d\Phi}{du}(u) \right],
\end{equation}
with 
$$
u=\frac{f_{\mathrm{min}}-\predmean(\x)}{\predvar(\x)}\:
$$
and $\Phi$ the normal cumulative distribution function.  The new
evaluation point is then chosen as a global maximizer of EI($\x$). An
example is given on Figure~\ref{fig:EI}, where the problem that
deceived the naive method of Figure~\ref{fig:optimsimple} is
directly solved with the EI criterion. First introduced in
(\cite{Jones1998}), this method has been used for computer experiments
in (\cite{Sasena2002}), while modified criteria have been used in
(\cite{Huang2005}) and (\cite{Williams2000}) to deal with noisy functions.

In (\cite{Jones2001}) and (\cite{Watson1995}), a fair number of
alternative criteria are presented and compared. Although quite
different in their formulation, they generally aim to answer the same
question: What is the most likely position of $\minimizer$?  Another,
and probably more relevant, question is: Where should the evaluation
be carried out optimally to improve our knowledge on the global
minimizers? 

In what follows, a criterion that addresses this question
will be presented, along with its performances. The reference for
comparison will be EI, which is, according to
\cite{Sasena2002}, a reasonable compromise between local and global
search, and has been successfully used in many applications.
\section{Estimating the density of  $\x^*$}
\label{sec:krig-cond-simul}
Once a Kriging surrogate model $\predmean$ has been obtained, any
global minimizer of $\predmean$ is a natural approximation of
$\minimizer$. However, it might be excessively daring to trust this
approximation as it does not take in account the uncertainty of the
prediction. A more cautious approach to estimating
$\minimizer$ is to use the probabilistic framework associated with
$F$. Of course, $\minimizer$ is not necessarily unique, and we shall
focus on describing the set of all global minimizers of $f$ as
efficiently as possible.
\subsection{Probabilistic modeling of the global minimizers of $f$}
\label{sec:minimizer-as-random-variable} 
According to the GP model, a global minimizer $\minimizer$ of $f$
corresponds to a global minimizer of this particular sample path of
$F$. More formally, consider the \emph{random} set
$\MinimizerSet_{\XX}$ of the global minimizers of $F$ over $\XX$, i.e.
the set of all global minimizers for each sample path, which can be written
as
$$
\MinimizerSet_{\XX}(\omega)=\{\minimizer \: \in \: \XX |
F(\omega,\minimizer)=\min_{\bm{u} \: \in \: \XX} F(\omega,u)\}, \: \forall \omega \in \Omega,
$$
with $\Omega$ the sample set. To ensure that $\MinimizerSet_{\XX}$ is not
empty, we assume that $F$ has continuous sample paths with probability
one.  This continuity can be ensured through a proper choice of
covariances (see, e.g., (\cite{Abrahamsen1997})).

Let $\Minimizer$ be a random vector uniformly distributed on
$\MinimizerSet_{\XX}$. The probability density function of this random
vector conditional to past evaluation results, that we shall
thereafter call conditional density of the global minimizers, is of
great interest, as it allows one not only to estimate the global
minimizers of $f$ (for example, through the maximization of their
probability density function), but also to characterize the
uncertainty associated with this estimation. In fact, given the
results $\bm{f}_\Sample$ of previous evaluations, the probability
density function $p_{\Minimizer|\bm{f}_\Sample}(\x)$ of $\Minimizer$
conditionally to $\bm{f}_\Sample$ contains all of what has been
assumed and learned about the system.  However, no tractable
analytical expression for such a quantity is available
(\cite{Adler2000,Sjo2000}). To overcome this
difficulty, the approach taken here is to consider a discrete version
of the conditional distribution, and to approximate it using Monte
Carlo simulations.

Let $\grid=\{\x_1, \ldots,\x_{\ngrid} \}$ be a finite subset of $\XX$,
$\MinimizerSet_{\grid}$ be the random set of global minimizers of $F$
over $\grid$, and $\Minimizergrid$ be a random vector uniformly
distributed on $\MinimizerSet_{\grid}$. The conditional probability mass function
of $\Minimizergrid$ given $\bm{f}_\Sample$ is
then $\forall \x \in \grid$
$$ 
\ddpcond(\x)\:=\:\Prob(\Minimizergrid=\x\,|\,\bm{F}_\Sample=\bm{f}_\Sample )\;.    
$$ 

It can be approximated using conditional simulations, \emph{i.e.},
simulations of $F$ that satisfy $\bm{F}_\Sample=\bm{f}_\Sample$.
Assuming that non-conditional simulations are available, several
methods exist to make them conditional (\cite{Chiles1999}).
Conditioning by Kriging seems the most promising of them in the
present context and will be presented in the next section.
\subsection{Conditioning by Kriging}
\label{sec:conditionning-by-Kriging} 
This method, due to G. Matheron, uses the unbiasedness of the Kriging
prediction to transform non-conditional simulations into simulations
interpolating the results $\bm{f}_{\Sample}$ of the evaluations. Let
$Z$ be a zero-mean Gaussian process with covariance $k$ (the same as
that of $F$) and $\hat{Z}$ be its Kriging predictor based on the random
variables $Z(\x_i),\:\x_i\in\Sample$, and consider the random process
\begin{equation}
  \label{eq:simucond}
   T(\x) =\predmean(\x)+ \left[Z(\x)-\hat{Z}(\x) \right],
\end{equation}
where $\predmean$ is the mean of the Kriging predictor based on the
design points in $\Sample$. Since this Kriging predictor is an
interpolator, at evaluation points in $\Sample$, we have
$\predmean(\x_i)=f(\x_i)$. Equation
(\ref{eq:prediction_at_evaluation_points}) implies that
$Z(\x_i)=\hat{Z}(\x_i)$, which leads to $T(\x_i)=f(\x_i),\: \forall \x_i \in
\Sample$. In other words, $T$ is such that all its sample paths
interpolate the known values of $f$.  It is then easy to
check that $T$ has the same finite-dimension distributions as $F$
conditionally to past evaluation results (\cite{Delfiner1977}).
Note that the same vector $\lambdas(\x)$ of Kriging coefficients is
used to interpolate the data and the simulations at design points.
Using (\ref{eq:Pred}), one can rewrite (\ref{eq:simucond}) as
\begin{equation}
  \label{eq:simucondtop}
  T(\x)=Z(\x)+\lambdas(\x)\Tr\left[ \bm{f}_\Sample- \bm{Z_{\Sample}}\right],
\end{equation}
with $\bm{Z_{\Sample}}=\lbrack Z(\x_1),\ldots,Z(\x_{\nobs}) \rbrack\Tr$.

In summary, to simulate $F$ over $\grid$ conditionally to past
evaluation results $\bm{f_\Sample}$, we can simulate a zero-mean
Gaussian process $Z$ over $\grid$, and use the following method:
\begin{itemize}
\item Compute, for every point in $\grid$, the vector of Kriging
  coefficients based on the design points in $\Sample$,
\item compute the Kriging prediction $\predmean(\x)$ based on past
  evaluation results $\bm{f}_\Sample$ for every $\x$ in $\grid$,
\item sample over $\grid$ non-conditional sample paths of $Z$
  (provided that a Gaussian sampler is available, setting the proper
  covariance can be achieved using, for example, the Cholesky
  decomposition),
\item
apply (\ref{eq:simucondtop}) at every point in $\grid$.
\end{itemize}

\begin{figure}
  \centering
  % This file is generated by the MATLAB m-file laprint.m. It can be included
% into LaTeX documents using the packages graphicx, color and psfrag.
% It is accompanied by a postscript file. A sample LaTeX file is:
%    \documentclass{article}\usepackage{graphicx,color,psfrag}
%    \begin{document}\input{simucond}\end{document}
% See http://www.mathworks.de/matlabcentral/fileexchange/loadFile.do?objectId=4638
% for recent versions of laprint.m.
%
% created by:           LaPrint version 3.16 (13.9.2004)
% created on:           20-Jun-2006 16:45:18
% eps bounding box:     15 cm x 10.7813 cm
% comment:              
%
\begin{psfrags}%
\psfragscanon%
%
% text strings:
\psfrag{s03}[b][b]{\color[rgb]{0,0,0}\setlength{\tabcolsep}{0pt}\begin{tabular}{c}$f$,$\predmean$\end{tabular}}%
\psfrag{s07}[b][b]{\color[rgb]{0,0,0}\setlength{\tabcolsep}{0pt}\begin{tabular}{c}$z$,$\hat{z}$\end{tabular}}%
\psfrag{s10}[t][t]{\color[rgb]{0,0,0}\setlength{\tabcolsep}{0pt}\begin{tabular}{c}$\x$\end{tabular}}%
\psfrag{s11}[b][b]{\color[rgb]{0,0,0}\setlength{\tabcolsep}{0pt}\begin{tabular}{c}$\predmean+(z-\hat{z})$,$\predmean$\end{tabular}}%
%
% xticklabels:
\psfrag{x01}[t][t]{0}%
\psfrag{x02}[t][t]{0.1}%
\psfrag{x03}[t][t]{0.2}%
\psfrag{x04}[t][t]{0.3}%
\psfrag{x05}[t][t]{0.4}%
\psfrag{x06}[t][t]{0.5}%
\psfrag{x07}[t][t]{0.6}%
\psfrag{x08}[t][t]{0.7}%
\psfrag{x09}[t][t]{0.8}%
\psfrag{x10}[t][t]{0.9}%
\psfrag{x11}[t][t]{1}%
\psfrag{x12}[t][t]{0}%
\psfrag{x13}[t][t]{0.1}%
\psfrag{x14}[t][t]{0.2}%
\psfrag{x15}[t][t]{0.3}%
\psfrag{x16}[t][t]{0.4}%
\psfrag{x17}[t][t]{0.5}%
\psfrag{x18}[t][t]{0.6}%
\psfrag{x19}[t][t]{0.7}%
\psfrag{x20}[t][t]{0.8}%
\psfrag{x21}[t][t]{0.9}%
\psfrag{x22}[t][t]{1}%
\psfrag{x23}[t][t]{0}%
\psfrag{x24}[t][t]{0.2}%
\psfrag{x25}[t][t]{0.4}%
\psfrag{x26}[t][t]{0.6}%
\psfrag{x27}[t][t]{0.8}%
\psfrag{x28}[t][t]{1}%
\psfrag{x29}[t][t]{0}%
\psfrag{x30}[t][t]{0.2}%
\psfrag{x31}[t][t]{0.4}%
\psfrag{x32}[t][t]{0.6}%
\psfrag{x33}[t][t]{0.8}%
\psfrag{x34}[t][t]{1}%
\psfrag{x35}[t][t]{0}%
\psfrag{x36}[t][t]{0.2}%
\psfrag{x37}[t][t]{0.4}%
\psfrag{x38}[t][t]{0.6}%
\psfrag{x39}[t][t]{0.8}%
\psfrag{x40}[t][t]{1}%
%
% yticklabels:
\psfrag{v01}[r][r]{0}%
\psfrag{v02}[r][r]{0.1}%
\psfrag{v03}[r][r]{0.2}%
\psfrag{v04}[r][r]{0.3}%
\psfrag{v05}[r][r]{0.4}%
\psfrag{v06}[r][r]{0.5}%
\psfrag{v07}[r][r]{0.6}%
\psfrag{v08}[r][r]{0.7}%
\psfrag{v09}[r][r]{0.8}%
\psfrag{v10}[r][r]{0.9}%
\psfrag{v11}[r][r]{1}%
\psfrag{v12}[r][r]{0}%
\psfrag{v13}[r][r]{0.1}%
\psfrag{v14}[r][r]{0.2}%
\psfrag{v15}[r][r]{0.3}%
\psfrag{v16}[r][r]{0.4}%
\psfrag{v17}[r][r]{0.5}%
\psfrag{v18}[r][r]{0.6}%
\psfrag{v19}[r][r]{0.7}%
\psfrag{v20}[r][r]{0.8}%
\psfrag{v21}[r][r]{0.9}%
\psfrag{v22}[r][r]{1}%
\psfrag{v23}[r][r]{9}%
\psfrag{v24}[r][r]{13}%
\psfrag{v25}[r][r]{-2}%
\psfrag{v26}[r][r]{2}%
\psfrag{v27}[r][r]{9}%
\psfrag{v28}[r][r]{13}%
%
% Figure:
\resizebox{12cm}{!}{\includegraphics{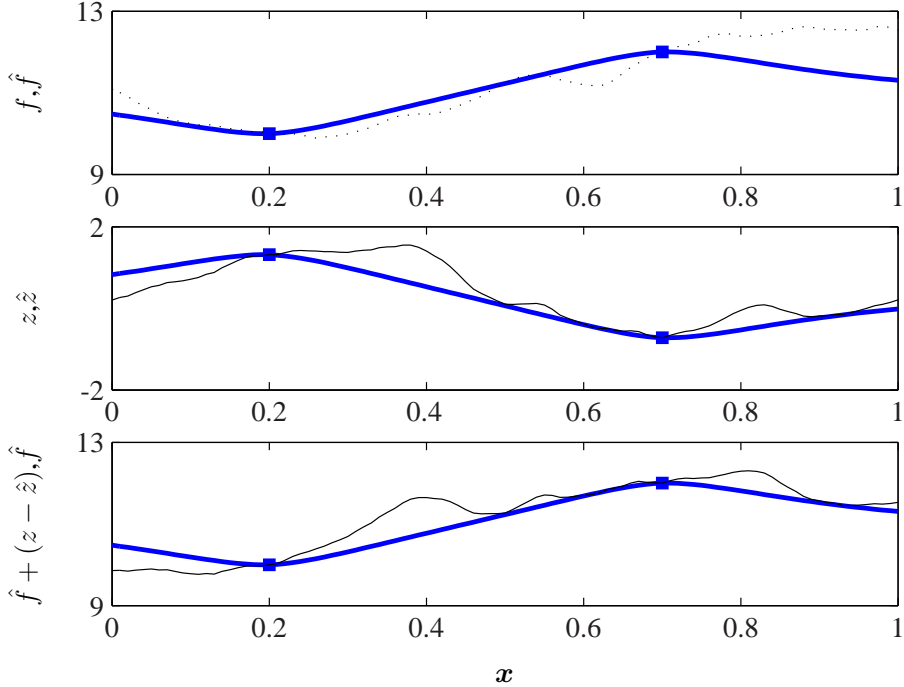}}%
\end{psfrags}%
%
% End simucond.tex

  \caption{Conditioning a simulation: (\emph{top}) unknown real curve $f$
    (doted line), sample points (squares) and associated Kriging
    prediction $\predmean$ (bold line); (\emph{middle}) non-conditional
    simulation $z$, sample points and associated Kriging prediction
    $\hat{z}$ (bold line); (\emph{bottom}) the simulation of the Kriging
    error $z-\hat{z}$ is picked up from the non-conditional simulation
    and added to the Kriging prediction to get the conditional
    simulation (thin line).}
  \label{fig:simucond}
\end{figure}
With this sampling method (see Figure~\ref{fig:simucond} for an
illustration), it becomes straightforward to estimate $\ddpcond$. Let
$\minimizer_i$ be a global minimizer of the $i$-th conditional
simulation ($i=1,\ldots,\ntraj$) over $\grid$ (if it is not unique, choose
one randomly).  Then, for any $\x$ in $\grid$ the empirical
probability mass distribution
\begin{equation}
  \label{eq:ddpapprox}
  \ddpaprox(\x)=\frac{1}{\ntraj}\sum_{i=1}^{\ntraj}
  \delta_{\minimizer_i}(\x), 
\end{equation}
with $\delta$ the Kronecker symbol, tends almost surely towards
$\ddpcond(\x)$ as $\ntraj$ tends towards infinity.  Moreover,
$\ddpcond$ tends in distribution towards
$p_{\Minimizer|\bm{f}_\Sample}$ as $\grid$ becomes dense in $\XX$.

\begin{figure}
  \centering
  % This file is generated by the MATLAB m-file laprint.m. It can be included
% into LaTeX documents using the packages graphicx, color and psfrag.
% It is accompanied by a postscript file. A sample LaTeX file is:
%    \documentclass{article}\usepackage{graphicx,color,psfrag}
%    \begin{document}\input{ddpmin}\end{document}
% See http://www.mathworks.de/matlabcentral/fileexchange/loadFile.do?objectId=4638
% for recent versions of laprint.m.
%
% created by:           LaPrint version 3.16 (13.9.2004)
% created on:           20-Jun-2006 16:48:33
% eps bounding box:     15 cm x 10.7813 cm
% comment:              
%
\begin{psfrags}%
\psfragscanon%
%
% text strings:
\psfrag{s01}[t][t]{\color[rgb]{0,0,0}\setlength{\tabcolsep}{0pt}\begin{tabular}{c}\end{tabular}}%
\psfrag{s02}[b][b]{\color[rgb]{0,0,0}\setlength{\tabcolsep}{0pt}\begin{tabular}{c}$\predmean$\end{tabular}}%
\psfrag{s05}[t][t]{\color[rgb]{0,0,0}\setlength{\tabcolsep}{0pt}\begin{tabular}{c}$\x$\end{tabular}}%
\psfrag{s06}[b][b]{\color[rgb]{0,0,0}\setlength{\tabcolsep}{0pt}\begin{tabular}{c}$\ddpaprox$\end{tabular}}%
%
% xticklabels:
\psfrag{x01}[t][t]{0}%
\psfrag{x02}[t][t]{0.1}%
\psfrag{x03}[t][t]{0.2}%
\psfrag{x04}[t][t]{0.3}%
\psfrag{x05}[t][t]{0.4}%
\psfrag{x06}[t][t]{0.5}%
\psfrag{x07}[t][t]{0.6}%
\psfrag{x08}[t][t]{0.7}%
\psfrag{x09}[t][t]{0.8}%
\psfrag{x10}[t][t]{0.9}%
\psfrag{x11}[t][t]{1}%
\psfrag{x12}[t][t]{0}%
\psfrag{x13}[t][t]{0.1}%
\psfrag{x14}[t][t]{0.2}%
\psfrag{x15}[t][t]{0.3}%
\psfrag{x16}[t][t]{0.4}%
\psfrag{x17}[t][t]{0.5}%
\psfrag{x18}[t][t]{0.6}%
\psfrag{x19}[t][t]{0.7}%
\psfrag{x20}[t][t]{0.8}%
\psfrag{x21}[t][t]{0.9}%
\psfrag{x22}[t][t]{1}%
\psfrag{x23}[t][t]{0}%
\psfrag{x24}[t][t]{1}%
\psfrag{x25}[t][t]{2}%
\psfrag{x26}[t][t]{3}%
\psfrag{x27}[t][t]{4}%
\psfrag{x28}[t][t]{5}%
\psfrag{x29}[t][t]{6}%
\psfrag{x30}[t][t]{7}%
\psfrag{x31}[t][t]{7.98}%
\psfrag{x32}[t][t]{0}%
\psfrag{x33}[t][t]{1}%
\psfrag{x34}[t][t]{2}%
\psfrag{x35}[t][t]{3}%
\psfrag{x36}[t][t]{4}%
\psfrag{x37}[t][t]{5}%
\psfrag{x38}[t][t]{6}%
\psfrag{x39}[t][t]{7}%
\psfrag{x40}[t][t]{8}%
%
% yticklabels:
\psfrag{v01}[r][r]{0}%
\psfrag{v02}[r][r]{0.1}%
\psfrag{v03}[r][r]{0.2}%
\psfrag{v04}[r][r]{0.3}%
\psfrag{v05}[r][r]{0.4}%
\psfrag{v06}[r][r]{0.5}%
\psfrag{v07}[r][r]{0.6}%
\psfrag{v08}[r][r]{0.7}%
\psfrag{v09}[r][r]{0.8}%
\psfrag{v10}[r][r]{0.9}%
\psfrag{v11}[r][r]{1}%
\psfrag{v12}[r][r]{0}%
\psfrag{v13}[r][r]{0.1}%
\psfrag{v14}[r][r]{0.2}%
\psfrag{v15}[r][r]{0.3}%
\psfrag{v16}[r][r]{0.4}%
\psfrag{v17}[r][r]{0.5}%
\psfrag{v18}[r][r]{0.6}%
\psfrag{v19}[r][r]{0.7}%
\psfrag{v20}[r][r]{0.8}%
\psfrag{v21}[r][r]{0.9}%
\psfrag{v22}[r][r]{1}%
\psfrag{v23}[r][r]{0}%
\psfrag{v24}[r][r]{0.004}%
\psfrag{v25}[r][r]{}%
\psfrag{v26}[r][r]{0.012}%
\psfrag{v27}[r][r]{-2}%
\psfrag{v28}[r][r]{4}%
\psfrag{v29}[r][r]{10}%
%
% Figure:
\resizebox{12cm}{!}{\includegraphics{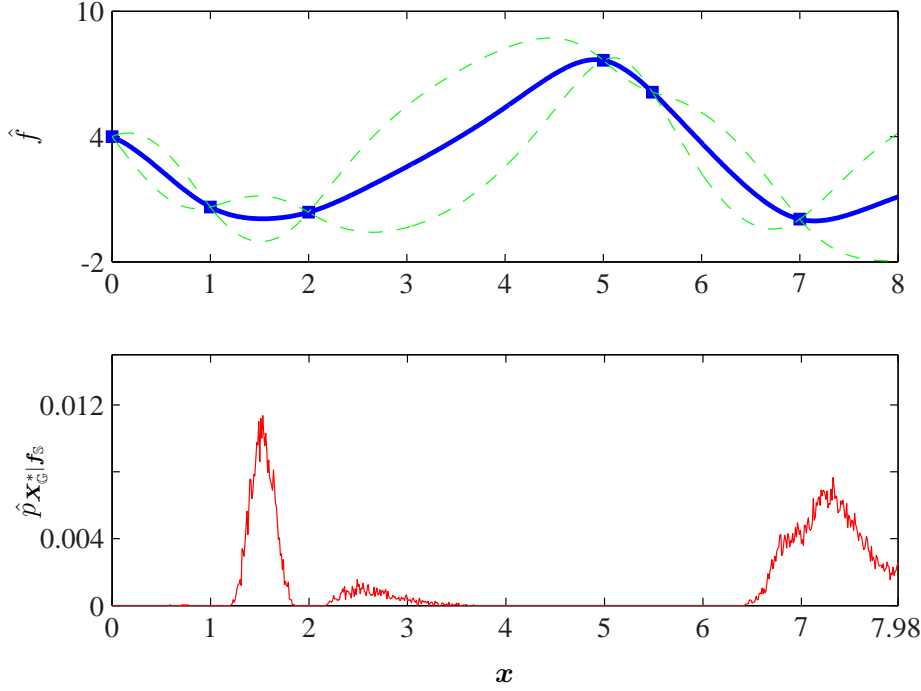}}%
\end{psfrags}%
%
% End ddpmin.tex

  \caption{Estimation of the conditional point mass distribution of
    $\Minimizer_{\grid}$: (\emph{top}) Kriging interpolation, 95\%
    confidence intervals and sample points; (\emph{bottom}) estimated
    point mass distribution of $\Minimizer_{\grid}$ using 10000 conditional
    simulations of $F$ and a regular grid for $\grid$.}
  \label{fig:ddpmin}
\end{figure}

Figure~\ref{fig:ddpmin} presents the approximation achieved by
$\ddpaprox$ for an example where locating a global minimizer is not
easy. Knowing the conditional distribution of $\Minimizer_\grid$ gives
valuable information on the areas of $\XX$ where a global minimizer
might be located, and that ought to be investigated.  This idea will
be detailed in the next section.

\section{The stepwise uncertainty reduction strategy}
\label{sec:stepw-uncert-reduct}
The knowledge about the global minimizers of $f$ is summarized by the
approximation $\ddpaprox$ of the conditional probability mass function
of $\Minimizergrid$. In order to evaluate the interest of a new
evaluation of $f$ at a given point, a measure of the expected
information gain is required.  An efficient measure is
\emph{conditional entropy}, as used in sequential testing
(\cite{Geaman1995}) in the \emph{Stepwise Uncertainty Reduction} (SUR)
strategy. This section extends the SUR strategy to global
optimization.

\subsection{Conditional entropy}
The entropy of a discrete random variable $U$ (expressed in bits) is
defined as:
$$
H(U)\:=\:-\,\sum_{u} \Prob(U=u)\log_2\, \Prob(U=u).
$$ 
$H(U)$ measures the spread of the distribution of $U$. It decreases
as this distribution gets more peaked. In particular :
\begin{itemize}
\item $\ddpaprox(\x)=1/\ngrid \;\; \forall \x \in \grid \,\Rightarrow
  \,H(\Minimizergrid)=\log_2(\ngrid)$,
\item
%  \begin{minipage}{1.0\linewidth}
$   
\ddpaprox(\x)=\left\{
  \begin{array}{lll}
    0 & \mbox{ if } & \x\:\neq\:\x_0 \\
    1 & \mbox{ if } & \x\:=\:\x_0
  \end{array} \right.\,\Rightarrow
  \,H(\Minimizergrid)=0
$
%  \end{minipage}
\end{itemize}
Similarly, for any event $\event$, the entropy of $U$ relative to the
probability measure $\Prob(.|\event)$ is 
$$
H(U|\event)\:=\:-\,\sum_{u} \Prob(U=u|\event) \log_2\, \Prob(U=u|\event).
$$
The conditional entropy of U given another discrete random variable $V$ is 
$$
H(U|V)\:=\:\sum_{v}\Prob(V=v)H(U|V=v),
$$
and the conditional entropy of $U$ given $\event$ and $V$ is
\begin{equation}
\label{eq:entropycond}
H(U|\event,V)\:=\:\sum_{v}\Prob(V=v|\event)H(U|\event,V=v).
\end{equation}
More details on conditional entropy can be found in \cite{Cover1991}.

\subsection{Conditional minimizer entropy}

Consider $F_Q(\x)$ a discrete version of $F(\x)$, defined as
$F_Q(\x)=Q(F(\x))$ with $Q$ a quantization operator. $Q$ is
characterized by a finite set of $M$ real numbers $\{y_1, \ldots,
y_{\npart}\}$, and defined as
$$
\forall u\:\in\: \mathbb{R}\:\: Q(u)=\left\{\begin{array}{l} y_1 \mbox{ if } u\leq
    y_1 \\ y_i \mbox{ if } y_{i-1}<u\leq y_i\;\;\forall i\:\in\: \llbracket2,\npart\rrbracket \\
    y_{\npart} \mbox{ if } y_{\npart}<u \end{array}\right.
$$
For optimization problems, the SUR strategy for the selection of the
next value of $\x \in \XX$ at which $f$ will be evaluated will be based
on $H(\Minimizergrid|\bm{F}_{\Sample}=\bm{f}_\Sample,F_Q(\x))$, the conditional entropy of $\Minimizergrid$ given $F_Q(\x)$ and
the evaluation results $\{\bm{F}_{\Sample}=\bm{f}_\Sample\}$ (we
shall refer to it later on as conditional entropy of the minimizer, or
simply conditional entropy).

Using (\ref{eq:entropycond}) we can write
\begin{equation}
  \label{eq:CEMcomputation}
  H(\Minimizergrid|\bm{F}_{\Sample}=\bm{f}_\Sample,F_Q(\x))\:=\:\sum_{i=1}^{\npart}\Prob(F_Q(\x)=y_i|\bm{F}_{\Sample}=\bm{f}_\Sample)H(\Minimizergrid|\bm{F}_{\Sample}=\bm{f}_\Sample,F_Q(\x)=y_i)
\end{equation}
with 
$$
H(\Minimizergrid|\bm{F}_{\Sample}=\bm{f}_\Sample,F_Q(\x)=y_i)\:=\:-\,\sum_{\bf{u}
  \in \grid} p_{\Minimizer_{\grid}|\bm{f}_\Sample,y_i}({\bf{u}}) \log_2\,
  p_{\Minimizer_{\grid}|\bm{f}_\Sample,y_i}({\bf{u}})\; ,
$$
and
$$p_{\Minimizer_{\grid}|\bm{f}_\Sample,y_i}(\bm{u})=\Prob(\Minimizer=\bm{u}
| \bm{F}_{\Sample}=\bm{f}_\Sample,F_Q(\x)=y_i ).$$ 

$H(\Minimizergrid|\bm{F}_{\Sample}=\bm{f}_\Sample,F_Q(\x))$ is a
measure of the anticipated uncertainty remaining in
$\Minimizergrid$ given the candidate evaluation point $\x$ and the
result $\bm{f}_\Sample$ of the previous evaluations. Anticipation is
introduced in (\ref{eq:CEMcomputation}) by considering the entropy of
$\Minimizer_{\grid}$ resulting from every possible sample value of
$F_Q(\x)$. At each stage of the iterative optimization, the SUR
strategy retains for the next evaluation a point that minimizes the
expected conditional minimizer entropy after the evaluation,
\emph{i.e.}, a point that maximizes the expected gain in information
about $\Minimizergrid$.

The conditional entropy of the minimizer thus takes in account the
conditional statistical properties of $F$ and particularly the
covariance of the model. There lies the interest of the SUR strategy
applied to global optimization.  It makes use of what has been
previously assumed and learned about $f$ to pick up the most
informative evaluation point. By contrast, the EI criterion (as most
standard criteria) depends only on the conditional mean and variance
of $F$ at the design point considered.
\section{Implementing the SUR strategy}
\label{sec:impl-sur-strat}
\subsection{IAGO algorithm}

Our algorithm is similar in spirit to a particular strategy for
Kriging-based optimization known as \emph{Efficient Global
  Optimization} (EGO) (\cite{Jones1998}). EGO starts with a small
initial design, estimates the parameters of the covariance of $F$ and
computes the Kriging model. Based on this model, an additional point
is selected in the design space to be the location of the next
evaluation of $f$ using the EI criterion.  The parameters of the
covariance are then re-estimated, the model re-computed, and the
process of choosing new points continues until the improvement
expected from sampling additional points has become sufficiently
small. The IAGO algorithm uses the same idea of iterative incorporation
of the obtained information to the prior on the function, but with a
different criterion.

The computation of the conditional entropy using (\ref{eq:CEMcomputation})
requires the choice of a quantization operator $Q$. We use the fact
that $F(\x)$ is conditionally Gaussian with known mean and variance
(obtained by Kriging), to select a set of possible values $\{y_1, \ldots,
y_{\npart}\}$, such that
$$P(F_Q(\x)=y_i|\bm{F}_{\Sample}=\bm{f}_\Sample)=\frac{1}{\npart}\;
\forall\:i\in\:\llbracket1:\npart\rrbracket\:.$$

By doing so, we choose a different quantization operator $Q$ for each
value of $\x$ to improve the precision with which the empirical mean
of the entropy reduction over possible evaluation results is computed.
This is simply an improved version of Monte Carlo integration. Here we
used a set of ten possible values ($M=10$).

For each of these possible values (or hypotheses $F(\x)=y_i$),
$\hat{p}_{\bm{f}_\Sample,y_i}$ is computed using conditional
simulations. The conditional entropy is then obtained using
(\ref{eq:CEMcomputation}). These operations are carried out on a
discrete set of candidate evaluation points (see
Section~\ref{sec:computational-burden} for some details on the choice
of this set), and a new evaluation of $f$ is finally performed at a
point that minimizes the conditional entropy. Next, as in the EGO algorithm, the
covariance parameters are re-estimated and the model re-computed.  The
procedure for the choice of an additional evaluation point is
described in Table~\ref{tab:algorithm}.

\begin{table}
  \centering
  \begin{algor}{}{
    \qinput{Set $\Sample=\{\x_1,\ldots,\x_{\nobs}\}$ of evaluation points and corresponding
      values $\bm{f}_{\Sample}$ of the function $f$}
    \qoutput{Additional evaluation point $\x_{\mathrm{new}}$}}
  Choose $\grid$, a discrete representation of $\XX$ \\
  Set covariance parameters either a priori or by
  maximum-likelihood estimation based on  $\bm{f}_{\Sample}$\\
  Compute $\ntraj$ non-conditional simulations over $\grid$\\
  Compute $\predmean(\x)$ and $\predvar(\x)$ over $\grid$ by Kriging
  from  $\bm{f}_{\Sample}$ \\
  \qwhile the set of candidate points has not been entirely explored\\
  \qdo Take an untried point $\x_{\mathrm{c}}$ in the set of candidate points\\
  Compute the parameters $\{y_1, \ldots, y_{\npart}\}$ of the quantization
  operator $Q$\\
  Compute $\bm{\Lambda}=\left[\lambdas(\x_1), \ldots,\lambdas(\x_{\ngrid})
  \right]$ the matrix of Kriging coefficients at every point in $
  \grid$  based on  evaluation points in $\Sample$ and $\x_{\mathrm{c}}$\\
  \qfor $i \qlet 1$ \qto $\npart$\\
  \qdo Construct conditional simulations using
  (\ref{eq:simucondtop}) and assuming that $f(\x_{\mathrm{c}})=y_i$\\
  Find a global minimizer $\minimizer_k$ of the $k$th conditional simulation over $\grid$ ($k=1,\ldots,\ntraj$) \\
  Estimate $p_{\Minimizer_{\grid}|\bm{f}_\Sample,y_i}$ over $\grid$ using (\ref{eq:ddpapprox})\\
  Compute $H(\Minimizergrid|\bm{F_{\Sample}}\,=\,\bm{f_{\Sample}},F_Q(\x_{\mathrm{c}})=y_i)$ \qrof\\
  Compute the conditional entropy given an evaluation at $\x_{\mathrm{c}}$ using
  (\ref{eq:CEMcomputation}) \qend\\
  Output $\x_{\mathrm{new}}$ that minimizes the conditional entropy over the set of
  candidate points
  \end{algor}
  \caption{Selection of a new evaluation point for $f$.}
  \label{tab:algorithm}
\end{table}

When the number of additional function evaluations is not specified
beforehand, we propose to use as a stopping criterion the conditional
probability that the global minimum of the GP model be no further
apart of the current minimum $f_{\mathrm{min}}^n$ of the Kriging
interpolation than a given tolerance threshold $\delta$. The algorithm then
stops when
$$
\Prob(F^*<f_{\mathrm{min}}^n+\delta|\bm{F}_{\Sample}=\bm{f}_\Sample)<P_{\mathrm{Stop}}\: ,
$$
with $F^*=\minim_{\x\:\in\:\grid}F(\x)$, and $P_{\mathrm{Stop}} \in [0,1]$ a critical
value to be chosen by the user. Proposed by Schonlau
(\cite{Schonlau1997}), this stopping criterion is well suited here,
since evaluating the repartition function of $f(\minimizer)$ does not
require any additional computation.

\subsection{Computational complexity}
\label{sec:computational-burden}

With the previous notation, $\nobs$ the number of evaluation points,
$\ntraj$ the number of conditional simulations, $\ngrid$ the number of
points in $\grid$ and $M$ number of possible results for an
evaluation, the computational complexity for the approximation of conditional
entropy (Steps 7 to 14 in Table~\ref{tab:algorithm}) is as follows:
\begin{itemize}
\item computing Kriging coefficients at every point in $\grid$ (Step
  8): $\GO(\nobs\ngrid$) (once the covariance matrix in
    (\ref{eq:Kr1}) has been factorized,
    Kriging at an untried point is simply in $\GO(\nobs)$),
\item constructing conditional simulations (Step 10):
  $\GO(\nobs\ntraj\ngrid)$ ($M$ is actually not involved
  since the main part of conditioning can be carried out outside the
  loop on the possible evaluation values),
\item locating the global minimizers for each simulation by exhaustive
  search (Step~11): $\GO(\ntraj\ngrid M)$,
\end{itemize}
Since all other operations are in $\GO(\ngrid)$ at most, evaluating
conditional entropy at any given point requires $\GO(\ngrid)$
operations.

To complete the description of an implementable algorithm, we must
specify a choice for $\grid$ and a policy for the minimization of
conditional entropy. What follows is just an example of a possible
strategy, and many variants could be considered.

The simplest choice for $\grid$ is a uniform grid on $\XX$.  However,
as the number of evaluations of $f$ increases, the spread of
$\ddpcond$ diminishes along with the precision for the computation of
the entropy. To keep a satisfactory precision over time, $\grid$ can
be a random sample of points in $\XX$, re-sampled after every
evaluation of $f$ with the density $\ddpaprox$.  Re-sampling makes it
possible to use a set $\grid$ with a smaller cardinal and to escape,
at least partly, the curse of dimensionality (to resample using
  $\ddpaprox$, any non-parametric density estimator could be used
  along with a sampling method such as Metropolis Hastings, see, e.g.,
  (\cite{Chib1995})).

Ideally, to choose an additional evaluation point for $f$ using IAGO,
conditional entropy should be minimized over $\XX$.  However, this of
course is in itself a global optimization problem, with many local
optima. It would be possible to design an ad-hoc optimization method
(as in (\cite{Jones2001})), but this perspective is not explored
here.  Instead, we evaluate the criterion extensively over a chosen
set of candidate points. Note that only the surrogate model is
involved at this stage, which makes the approach practical. The
  idea is, exactly as for the choice of $\grid$, to use a
  space-filling sample covering $\XX$ and resampled after each new
  evaluation. The current implementation of IAGO simply uses a Latin
  Hyper Cube (LHC) sample, however, it would be easy to adapt this
  sample iteratively using the conditional distribution of the
  minimizers $\ddpaprox$ as a prior.  For instance, areas of the
  design space where the density is sufficiently small could be
  ignored.  After a few evaluations, a large portion of the design
  space usually satisfies this property, and the computations saved
  could be used to improve knowledge on the criterion by sampling
  where $\ddpaprox$ is high (using the same approach as for the choice
  of $\grid$).

As dimension increases, trying to cover the factor space while keeping
the same accuracy leads to an exponential increase in complexity.
However, in a context of expensive function evaluation, the objective
is less to specify exactly all global minimizers (which could be too
demanding in function evaluations anyway), than to use available
information efficiently to reduce the likely areas for the location of
these minimizers. This is exactly the driving concept behind IAGO.
In practice, within a set of one thousand candidate points, picking an
additional evaluation point requires about five minutes with a
standard personal computer (and this figure is relatively independent
of the dimension of the factor space). Moreover, the result obtained
can be trusted to be a consistent choice within this set of candidate
points, in regard of what has been assumed and learned about $f$.
\subsection{Taking noise in account}
\label{sec:taking-noise-account}
Practical optimization problems often involve noise. This
section discusses possible adaptations of the optimization
algorithm that make it possible to deal with noisy situations, namely
noise on the evaluation of $f$ and noise on the factors.
\subsubsection{Noise on the evaluation of $f$}
\label{sec:noise}
When the evaluations of $f$ are corrupted by noise, the algorithm must
take this fact into account. A useful tool to deal with such
situations is \emph{non-interpolative Kriging} (see
Section~\ref{sec:dealing-with-noisy}).

\begin{figure}
  \centering
  % This file is generated by the MATLAB m-file laprint.m. It can be included
% into LaTeX documents using the packages graphicx, color and psfrag.
% It is accompanied by a postscript file. A sample LaTeX file is:
%    \documentclass{article}\usepackage{graphicx,color,psfrag}
%    \begin{document}\input{ExempleKrN}\end{document}
% See http://www.mathworks.de/matlabcentral/fileexchange/loadFile.do?objectId=4638
% for recent versions of laprint.m.
%
% created by:           LaPrint version 3.16 (13.9.2004)
% created on:           01-Jun-2006 14:48:25
% eps bounding box:     15 cm x 10.7813 cm
% comment:              
%
\begin{psfrags}%
\psfragscanon%
%
% text strings:
\psfrag{s04}[t][t]{\color[rgb]{0,0,0}\setlength{\tabcolsep}{0pt}\begin{tabular}{c}\end{tabular}}%
%
% xticklabels:
\psfrag{x01}[t][t]{0}%
\psfrag{x02}[t][t]{0.1}%
\psfrag{x03}[t][t]{0.2}%
\psfrag{x04}[t][t]{0.3}%
\psfrag{x05}[t][t]{0.4}%
\psfrag{x06}[t][t]{0.5}%
\psfrag{x07}[t][t]{0.6}%
\psfrag{x08}[t][t]{0.7}%
\psfrag{x09}[t][t]{0.8}%
\psfrag{x10}[t][t]{0.9}%
\psfrag{x11}[t][t]{1}%
\psfrag{x12}[t][t]{0}%
\psfrag{x13}[t][t]{0.1}%
\psfrag{x14}[t][t]{0.2}%
\psfrag{x15}[t][t]{0.3}%
\psfrag{x16}[t][t]{0.4}%
\psfrag{x17}[t][t]{0.5}%
\psfrag{x18}[t][t]{0.6}%
\psfrag{x19}[t][t]{0.7}%
\psfrag{x20}[t][t]{0.8}%
\psfrag{x21}[t][t]{0.9}%
\psfrag{x22}[t][t]{1}%
%
% yticklabels:
\psfrag{v01}[r][r]{0}%
\psfrag{v02}[r][r]{0.1}%
\psfrag{v03}[r][r]{0.2}%
\psfrag{v04}[r][r]{0.3}%
\psfrag{v05}[r][r]{0.4}%
\psfrag{v06}[r][r]{0.5}%
\psfrag{v07}[r][r]{0.6}%
\psfrag{v08}[r][r]{0.7}%
\psfrag{v09}[r][r]{0.8}%
\psfrag{v10}[r][r]{0.9}%
\psfrag{v11}[r][r]{1}%
\psfrag{v12}[r][r]{0}%
\psfrag{v13}[r][r]{0.5}%
\psfrag{v14}[r][r]{1}%
\psfrag{v15}[r][r]{1.5}%
\psfrag{v16}[r][r]{2}%
\psfrag{v17}[r][r]{2.5}%
\psfrag{v18}[r][r]{3}%
%
% Figure:
\resizebox{12cm}{!}{\includegraphics{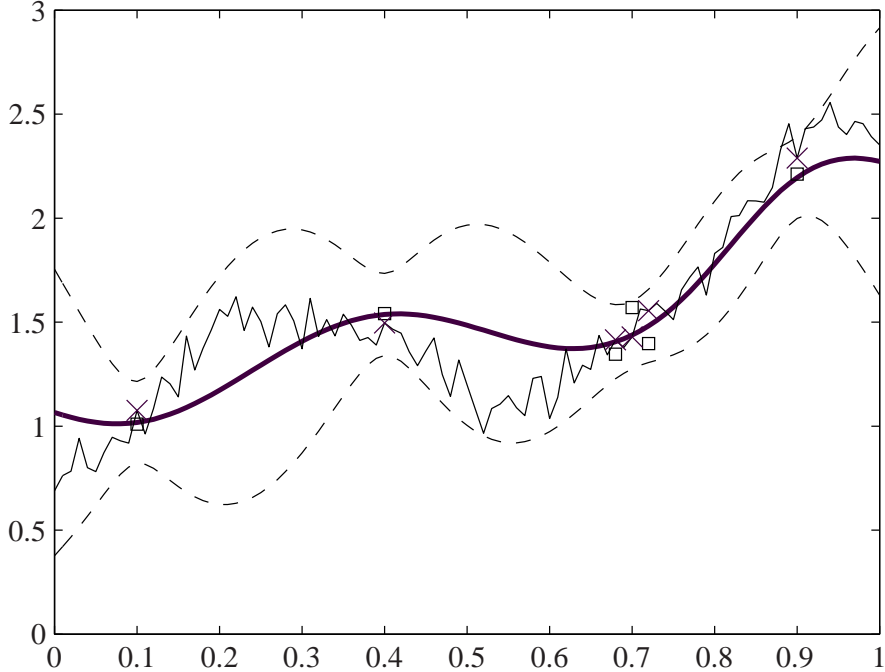}}%
\end{psfrags}%
%
% End ExempleKrN.tex

  \caption{Example of prediction by Kriging (bold line) of noisy
    measurements represented by squares. Dashed lines represent 95\%
    confidence regions for the prediction and the thin solid line is
    an example of conditional simulation obtained using the method
    presented in Section~\ref{sec:noise} (crosses represent the
    simulated measurements.)}
  \label{fig:ExempleKrN}
\end{figure}

If the evaluation at $\x_i \: \in\:\Sample$ is assumed to be
  corrupted by an additive Gaussian noise $\varepsilon_i$ with known mean and variance,
  the Kriging prediction should no longer be interpolative (see
  Section~\ref{sec:dealing-with-noisy}). The optimization algorithm
  remains nearly unchanged, except for the conditional simulations. We
  now wish to build sample paths of $F$ conditionally to evaluation
  results, i.e. realizations of the random variables
  $f(\x_i)+\varepsilon_i\mbox{ for } \x_i \in \Sample$. Since the variance of the
  prediction error is no longer zero at evaluation points (in other
  words, there is some uncertainty left on the values of $f$ at
  evaluation points), we first have to sample, at each evaluation
  points, from the distribution of $F$ conditionally to noisy
  evaluation results. An interpolative simulation, based on these
  sample values, is then built using conditioning by Kriging. An
  example of such a simulation is proposed on
  Figure~\ref{fig:ExempleKrN}.

\subsubsection{Noise on the factors}
In many industrial design problems, the variability of the values of
the factors in mass production has a significant impact on the
performance that can be achieved. In such a case, one might want to
design a system that optimizes some performance measure while ensuring
that the performance uncertainty (stemming from noise on the
  factors) remains under control. These so-called \emph{robust
  optimization} problems can generally be written as
\begin{equation}
\argmin_{\x\;\in\;\DD} J(\x)\:,
\label{eq:optimrobuste}
\end{equation}
with $J(\x)$ a cost function reflecting some statistical property of
the corrupted performance measure $f(\x+\bm{\varepsilon})$, where $\bm{\varepsilon}$ is a
random vector accounting for noise on the factors. Classical cost functions are:
%(ref):
\begin{itemize}
\item mean: $J(\x)=\expect_{\bm{\varepsilon}}[f(\x+\bm{\varepsilon})]$,
\item standard deviation: $J(\x)=\sigma_{\bm{\varepsilon}}[f(\x+\bm{\varepsilon})]=\sqrt(\var(f(\x+\bm{\varepsilon})))$,
\item linear combination of mean and standard deviation:
  $J(\x)=\expect_{\bm{\varepsilon}}[f(\x+\bm{\varepsilon})]+k\sigma_{\bm{\varepsilon}}[f(\x+\bm{\varepsilon})]$,
\item $\alpha$-quantile: $J(\x)=\quant(\x)$ \\ with $\quant(\x)$ such that
  $\Prob(f(\x+\bm{\varepsilon})<\quant(\x))=\alpha$.
\end{itemize}
Using, for example, the $\alpha$-quantile as a cost function, it is
possible to adapt our optimization algorithm to solve
(\ref{eq:optimrobuste}). Given a set of evaluation results $\bm{f_\Sample}$
at noise-free evaluation points, and assuming that it is possible to
sample from the distribution $p_{\bm{\varepsilon}}$ of $\bm{\varepsilon}$, a Monte-Carlo
approximation $\hat{Q}^{\alpha}(\x)$ of $\quant(\x)$ is easily obtained by
computing $\hat{f}(\x+\bm{\varepsilon})$ over a set sampled from $p_{\bm{\varepsilon}}$.
The global optimization algorithm can then be applied to $\quant(\x)$
instead of $f$, using pseudo-evaluations
$\bm{\hat{Q}^{\alpha}_{\Sample}}=[\hat{Q}^{\alpha}(\x_1),
\ldots,\hat{Q}^{\alpha}(\x_\nobs)]$ instead of $\bm{f_{\Sample}}$.

It is of course possible to combine these ideas and to deal
simultaneously with noise both on the factors and the function
evaluations.

\section{Illustrations}
\label{sec:illustrations}
This section contains some simple examples of global optimization
using IAGO, with a regular grid as a set of candidate evaluation
points. An empirical comparison with global optimization using
expected improvement is also presented. The Mat\'ern covariance class
will be used for Kriging prediction, as it facilitates the tuning of
the variance, regularity and range of correlation of the underlying
random process, but note that any kind of admissible covariance
  could have been used. The parameters of the covariance may be
estimated from the data using a maximum-likelihood approach (see
Section~\ref{sec:implementation}).

\subsection{A one-dimensional example}
\begin{figure}
    \centering \subfigure[Kriging prediction and point mass
    distribution of the global minimizers based on the initial
    design]{ % This file is generated by the MATLAB m-file laprint.m. It can be included
% into LaTeX documents using the packages graphicx, color and psfrag.
% It is accompanied by a postscript file. A sample LaTeX file is:
%    \documentclass{article}\usepackage{graphicx,color,psfrag}
%    \begin{document}\input{SUR1Dinitial}\end{document}
% See http://www.mathworks.de/matlabcentral/fileexchange/loadFile.do?objectId=4638
% for recent versions of laprint.m.
%
% created by:           LaPrint version 3.16 (13.9.2004)
% created on:           08-Jun-2006 14:58:16
% eps bounding box:     15 cm x 10.7813 cm
% comment:              
%
\begin{psfrags}%
\psfragscanon%
%
% text strings:
\psfrag{s06}[][]{\color[rgb]{0,0,0}\setlength{\tabcolsep}{0pt}\begin{tabular}{c} \end{tabular}}%
\psfrag{s07}[][]{\color[rgb]{0,0,0}\setlength{\tabcolsep}{0pt}\begin{tabular}{c} \end{tabular}}%
%
% xticklabels:
\psfrag{x01}[t][t]{0}%
\psfrag{x02}[t][t]{0.1}%
\psfrag{x03}[t][t]{0.2}%
\psfrag{x04}[t][t]{0.3}%
\psfrag{x05}[t][t]{0.4}%
\psfrag{x06}[t][t]{0.5}%
\psfrag{x07}[t][t]{0.6}%
\psfrag{x08}[t][t]{0.7}%
\psfrag{x09}[t][t]{0.8}%
\psfrag{x10}[t][t]{0.9}%
\psfrag{x11}[t][t]{1}%
\psfrag{x12}[t][t]{0}%
\psfrag{x13}[t][t]{2}%
\psfrag{x14}[t][t]{4}%
\psfrag{x15}[t][t]{6}%
%
% yticklabels:
\psfrag{v01}[r][r]{0}%
\psfrag{v02}[r][r]{0.1}%
\psfrag{v03}[r][r]{0.2}%
\psfrag{v04}[r][r]{0.3}%
\psfrag{v05}[r][r]{0.4}%
\psfrag{v06}[r][r]{0.5}%
\psfrag{v07}[r][r]{0.6}%
\psfrag{v08}[r][r]{0.7}%
\psfrag{v09}[r][r]{0.8}%
\psfrag{v10}[r][r]{0.9}%
\psfrag{v11}[r][r]{1}%
\psfrag{v12}[r][r]{-1}%
\psfrag{v13}[r][r]{1}%
\psfrag{v14}[r][r]{3}%
\psfrag{v15}[r][r]{5}%
\psfrag{v16}[r][r]{7}%
\psfrag{v17}[r][r]{9}%
%
% Figure:
\resizebox{6cm}{!}{\includegraphics{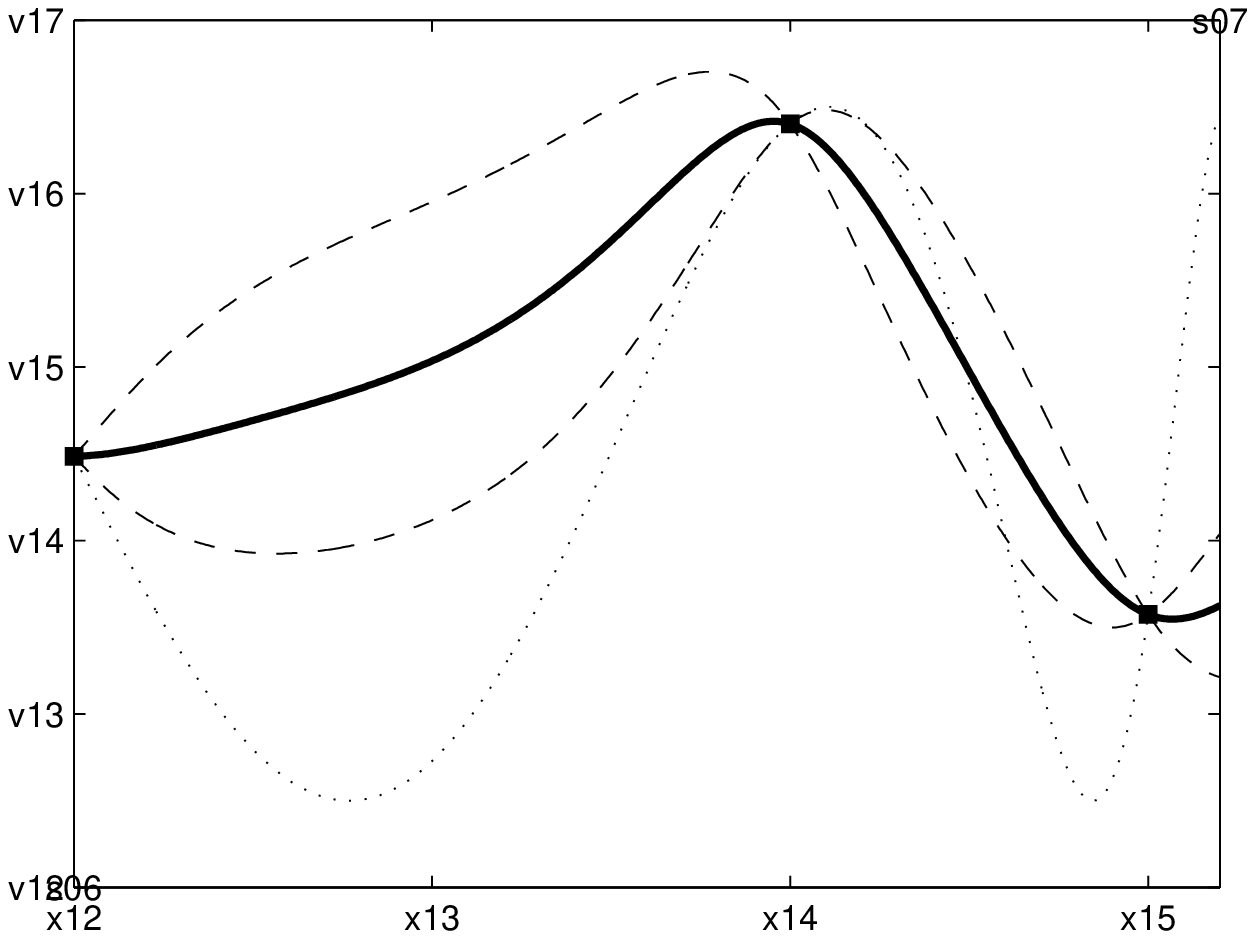}}%
\end{psfrags}%
%
% End SUR1Dinitial.tex

      \input{./courbes/SUR1Dinitialddp.tex}} \subfigure[Standard IAGO
    algorithm (noise free case)]{ % This file is generated by the MATLAB m-file laprint.m. It can be included
% into LaTeX documents using the packages graphicx, color and psfrag.
% It is accompanied by a postscript file. A sample LaTeX file is:
%    \documentclass{article}\usepackage{graphicx,color,psfrag}
%    \begin{document}\input{SUR1D}\end{document}
% See http://www.mathworks.de/matlabcentral/fileexchange/loadFile.do?objectId=4638
% for recent versions of laprint.m.
%
% created by:           LaPrint version 3.16 (13.9.2004)
% created on:           18-May-2006 09:51:20
% eps bounding box:     15 cm x 10.7813 cm
% comment:              
%
\begin{psfrags}%
\psfragscanon%
%
% text strings:
\psfrag{s01}[l][l]{\color[rgb]{0,0,0}\setlength{\tabcolsep}{0pt}\begin{tabular}{l}\large1\end{tabular}}%
\psfrag{s02}[l][l]{\color[rgb]{0,0,0}\setlength{\tabcolsep}{0pt}\begin{tabular}{l}\large2\end{tabular}}%
\psfrag{s03}[l][l]{\color[rgb]{0,0,0}\setlength{\tabcolsep}{0pt}\begin{tabular}{l}\large3\end{tabular}}%
\psfrag{s04}[l][l]{\color[rgb]{0,0,0}\setlength{\tabcolsep}{0pt}\begin{tabular}{l}\large4\end{tabular}}%
\psfrag{s05}[l][l]{\color[rgb]{0,0,0}\setlength{\tabcolsep}{0pt}\begin{tabular}{l}\large5\end{tabular}}%
\psfrag{s06}[l][l]{\color[rgb]{0,0,0}\setlength{\tabcolsep}{0pt}\begin{tabular}{l}\large6\end{tabular}}%
%
% xticklabels:
\psfrag{x01}[t][t]{0}%
\psfrag{x02}[t][t]{0.1}%
\psfrag{x03}[t][t]{0.2}%
\psfrag{x04}[t][t]{0.3}%
\psfrag{x05}[t][t]{0.4}%
\psfrag{x06}[t][t]{0.5}%
\psfrag{x07}[t][t]{0.6}%
\psfrag{x08}[t][t]{0.7}%
\psfrag{x09}[t][t]{0.8}%
\psfrag{x10}[t][t]{0.9}%
\psfrag{x11}[t][t]{1}%
\psfrag{x12}[t][t]{0}%
\psfrag{x13}[t][t]{2}%
\psfrag{x14}[t][t]{4}%
\psfrag{x15}[t][t]{6}%
%
% yticklabels:
\psfrag{v01}[r][r]{0}%
\psfrag{v02}[r][r]{0.1}%
\psfrag{v03}[r][r]{0.2}%
\psfrag{v04}[r][r]{0.3}%
\psfrag{v05}[r][r]{0.4}%
\psfrag{v06}[r][r]{0.5}%
\psfrag{v07}[r][r]{0.6}%
\psfrag{v08}[r][r]{0.7}%
\psfrag{v09}[r][r]{0.8}%
\psfrag{v10}[r][r]{0.9}%
\psfrag{v11}[r][r]{1}%
\psfrag{v12}[r][r]{-1}%
\psfrag{v13}[r][r]{1}%
\psfrag{v14}[r][r]{3}%
\psfrag{v15}[r][r]{5}%
\psfrag{v16}[r][r]{7}%
\psfrag{v17}[r][r]{9}%
%
% Figure:
\resizebox{6cm}{!}{\includegraphics{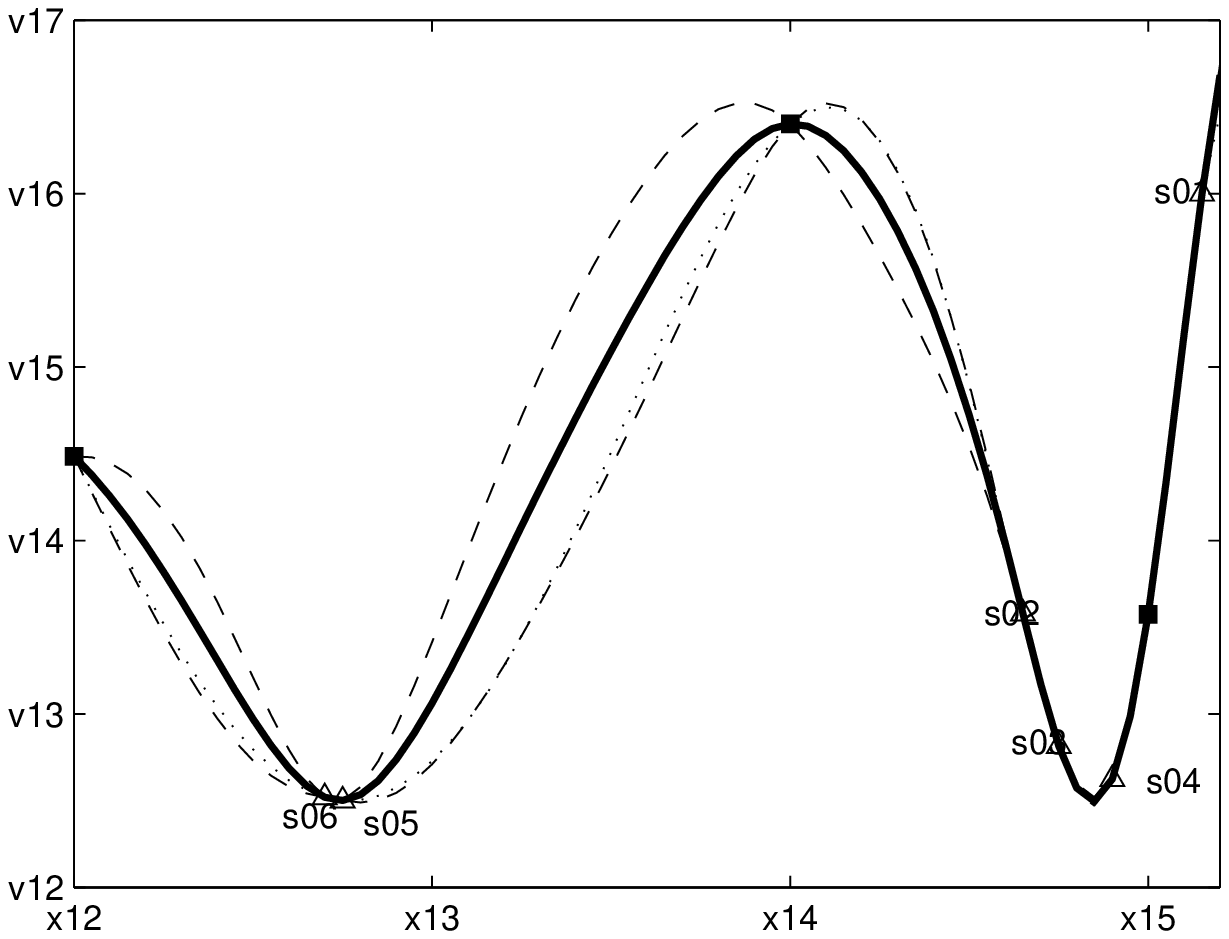}}%
\end{psfrags}%
%
% End SUR1D.tex

      \input{./courbes/SUR1Dddp.tex}} \subfigure[IAGO algorithm for
    noisy evaluations (the additive noise is zero-mean Gaussian with
    standard deviation 0.2)]{ % This file is generated by the MATLAB m-file laprint.m. It can be included
% into LaTeX documents using the packages graphicx, color and psfrag.
% It is accompanied by a postscript file. A sample LaTeX file is:
%    \documentclass{article}\usepackage{graphicx,color,psfrag}
%    \begin{document}\input{SUR1DNO}\end{document}
% See http://www.mathworks.de/matlabcentral/fileexchange/loadFile.do?objectId=4638
% for recent versions of laprint.m.
%
% created by:           LaPrint version 3.16 (13.9.2004)
% created on:           18-May-2006 10:01:09
% eps bounding box:     15 cm x 10.7813 cm
% comment:              
%
\begin{psfrags}%
\psfragscanon%
%
% text strings:
\psfrag{s01}[l][l]{\color[rgb]{0,0,0}\setlength{\tabcolsep}{0pt}\begin{tabular}{l}\large1\end{tabular}}%
\psfrag{s02}[l][l]{\color[rgb]{0,0,0}\setlength{\tabcolsep}{0pt}\begin{tabular}{l}\large2\end{tabular}}%
\psfrag{s03}[l][l]{\color[rgb]{0,0,0}\setlength{\tabcolsep}{0pt}\begin{tabular}{l}\large3\end{tabular}}%
\psfrag{s04}[l][l]{\color[rgb]{0,0,0}\setlength{\tabcolsep}{0pt}\begin{tabular}{l}\large4\end{tabular}}%
\psfrag{s05}[l][l]{\color[rgb]{0,0,0}\setlength{\tabcolsep}{0pt}\begin{tabular}{l}\large5\end{tabular}}%
\psfrag{s06}[l][l]{\color[rgb]{0,0,0}\setlength{\tabcolsep}{0pt}\begin{tabular}{l}\large6\end{tabular}}%
%
% xticklabels:
\psfrag{x01}[t][t]{0}%
\psfrag{x02}[t][t]{0.1}%
\psfrag{x03}[t][t]{0.2}%
\psfrag{x04}[t][t]{0.3}%
\psfrag{x05}[t][t]{0.4}%
\psfrag{x06}[t][t]{0.5}%
\psfrag{x07}[t][t]{0.6}%
\psfrag{x08}[t][t]{0.7}%
\psfrag{x09}[t][t]{0.8}%
\psfrag{x10}[t][t]{0.9}%
\psfrag{x11}[t][t]{1}%
\psfrag{x12}[t][t]{0}%
\psfrag{x13}[t][t]{2}%
\psfrag{x14}[t][t]{4}%
\psfrag{x15}[t][t]{6}%
%
% yticklabels:
\psfrag{v01}[r][r]{0}%
\psfrag{v02}[r][r]{0.1}%
\psfrag{v03}[r][r]{0.2}%
\psfrag{v04}[r][r]{0.3}%
\psfrag{v05}[r][r]{0.4}%
\psfrag{v06}[r][r]{0.5}%
\psfrag{v07}[r][r]{0.6}%
\psfrag{v08}[r][r]{0.7}%
\psfrag{v09}[r][r]{0.8}%
\psfrag{v10}[r][r]{0.9}%
\psfrag{v11}[r][r]{1}%
\psfrag{v12}[r][r]{-1}%
\psfrag{v13}[r][r]{1}%
\psfrag{v14}[r][r]{3}%
\psfrag{v15}[r][r]{5}%
\psfrag{v16}[r][r]{7}%
\psfrag{v17}[r][r]{9}%
%
% Figure:
\resizebox{6cm}{!}{\includegraphics{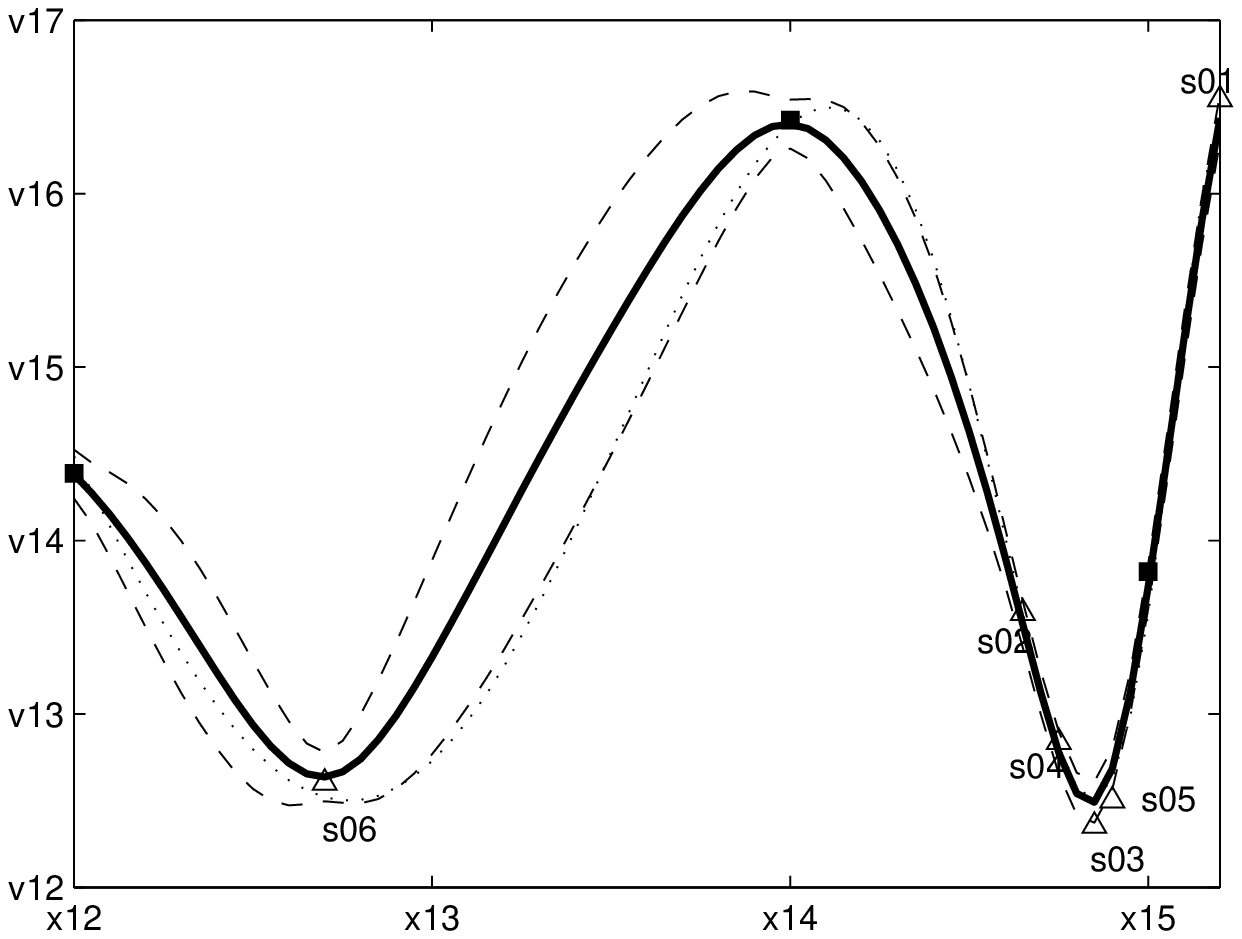}}%
\end{psfrags}%
%
% End SUR1DNO.tex

      \input{./courbes/SUR1DNOddp.tex}}
    \caption{Example of global optimization using IAGO on a function of
      one variable (dotted line), with an initial design consisting of
      three points (represented by squares).  Six additional
      evaluations are carried out (triangles) using two versions
      of the IAGO algorithm.  The graphs on the \emph{left part} of the
      figure account for final predictions, while the \emph{right
        part} presents the final point mass distributions of the
      global minimizers}
    \label{fig:SUR1D}
  \end{figure}

Consider the function with two global minimizers illustrated by Figure
\ref{fig:SUR1D} and defined by $f:
x\longmapsto4[1-\sin(x+8\exp(x-7))]$. Given an initial design
consisting of three points, the IAGO algorithm is used to compute six
additional points iteratively.  The final Kriging model is depicted in
the left part of Figure \ref{fig:SUR1D}, along with the resulting
point mass conditional distribution for the minimizer on the right
part. After adding some noise on the function evaluations, the
variant of the algorithm presented in Section~\ref{sec:noise} is
also applied to the function with the same initial design. In both
cases, the six additional evaluations have significantly reduced the
uncertainty associated with the position of the global minimizers. The
remaining likely locations reduce to two small areas centered on the
two actual global minimizers. In the noisy case, larger zones are
identified, a direct consequence of the uncertainty associated with
the evaluations.

\begin{figure}
    \centering
      % This file is generated by the MATLAB m-file laprint.m. It can be included
% into LaTeX documents using the packages graphicx, color and psfrag.
% It is accompanied by a postscript file. A sample LaTeX file is:
%    \documentclass{article}\usepackage{graphicx,color,psfrag}
%    \begin{document}\input{SUR1DNI}\end{document}
% See http://www.mathworks.de/matlabcentral/fileexchange/loadFile.do?objectId=4638
% for recent versions of laprint.m.
%
% created by:           LaPrint version 3.16 (13.9.2004)
% created on:           18-May-2006 10:14:14
% eps bounding box:     15 cm x 10.7813 cm
% comment:              
%
\begin{psfrags}%
\psfragscanon%
%
% text strings:
\psfrag{s01}[l][l]{\color[rgb]{0,0,0}\setlength{\tabcolsep}{0pt}\begin{tabular}{l}\large1\end{tabular}}%
\psfrag{s02}[l][l]{\color[rgb]{0,0,0}\setlength{\tabcolsep}{0pt}\begin{tabular}{l}\large2\end{tabular}}%
\psfrag{s03}[l][l]{\color[rgb]{0,0,0}\setlength{\tabcolsep}{0pt}\begin{tabular}{l}\large3\end{tabular}}%
\psfrag{s04}[l][l]{\color[rgb]{0,0,0}\setlength{\tabcolsep}{0pt}\begin{tabular}{l}\large4\end{tabular}}%
\psfrag{s05}[l][l]{\color[rgb]{0,0,0}\setlength{\tabcolsep}{0pt}\begin{tabular}{l}\large5\end{tabular}}%
\psfrag{s06}[l][l]{\color[rgb]{0,0,0}\setlength{\tabcolsep}{0pt}\begin{tabular}{l}\large6\end{tabular}}%
%
% xticklabels:
\psfrag{x01}[t][t]{0}%
\psfrag{x02}[t][t]{0.1}%
\psfrag{x03}[t][t]{0.2}%
\psfrag{x04}[t][t]{0.3}%
\psfrag{x05}[t][t]{0.4}%
\psfrag{x06}[t][t]{0.5}%
\psfrag{x07}[t][t]{0.6}%
\psfrag{x08}[t][t]{0.7}%
\psfrag{x09}[t][t]{0.8}%
\psfrag{x10}[t][t]{0.9}%
\psfrag{x11}[t][t]{1}%
\psfrag{x12}[t][t]{0}%
\psfrag{x13}[t][t]{2}%
\psfrag{x14}[t][t]{4}%
\psfrag{x15}[t][t]{6}%
%
% yticklabels:
\psfrag{v01}[r][r]{0}%
\psfrag{v02}[r][r]{0.1}%
\psfrag{v03}[r][r]{0.2}%
\psfrag{v04}[r][r]{0.3}%
\psfrag{v05}[r][r]{0.4}%
\psfrag{v06}[r][r]{0.5}%
\psfrag{v07}[r][r]{0.6}%
\psfrag{v08}[r][r]{0.7}%
\psfrag{v09}[r][r]{0.8}%
\psfrag{v10}[r][r]{0.9}%
\psfrag{v11}[r][r]{1}%
\psfrag{v12}[r][r]{-1}%
\psfrag{v13}[r][r]{1}%
\psfrag{v14}[r][r]{3}%
\psfrag{v15}[r][r]{5}%
\psfrag{v16}[r][r]{7}%
\psfrag{v17}[r][r]{9}%
%
% Figure:
\resizebox{6cm}{!}{\includegraphics{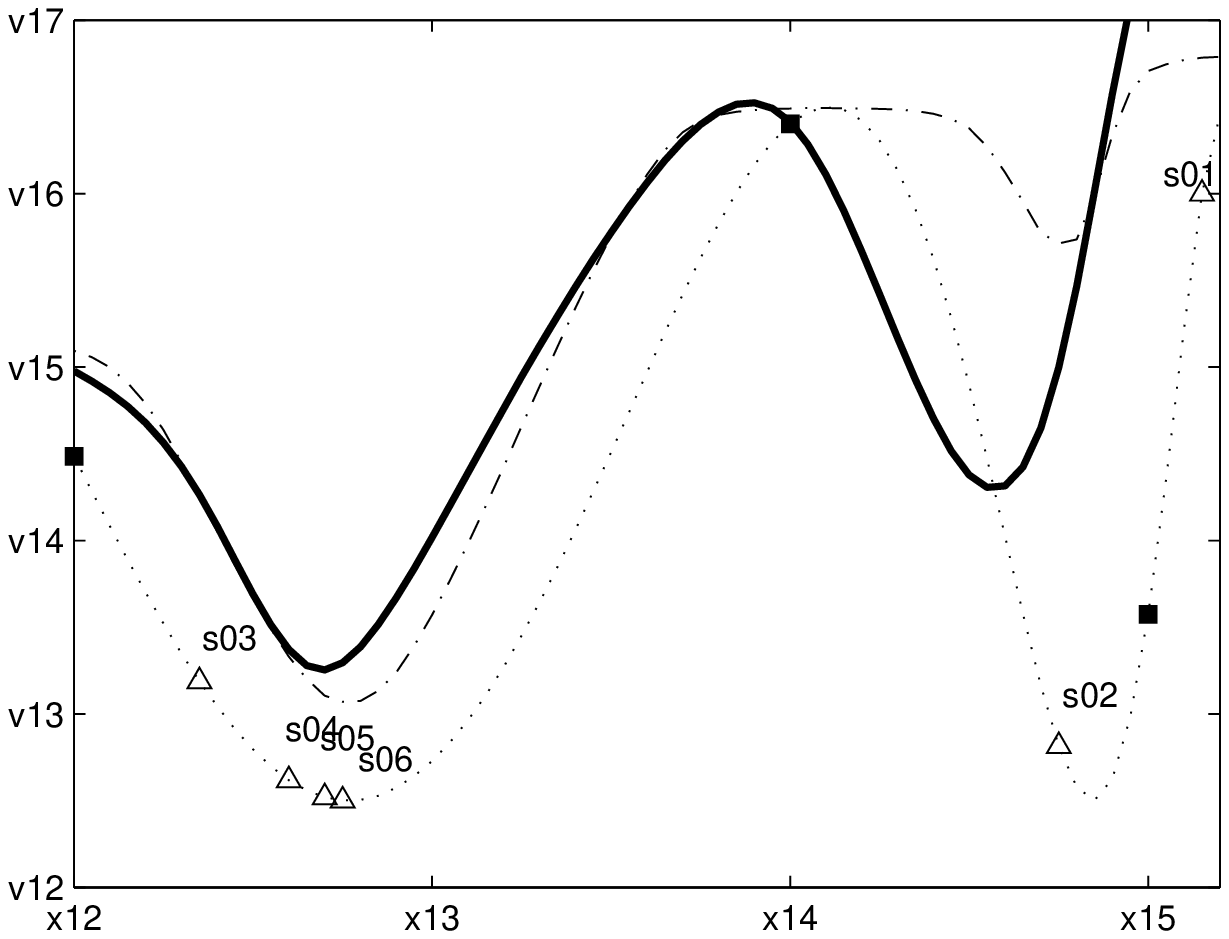}}%
\end{psfrags}%
%
% End SUR1DNI.tex

      \input{./courbes/SUR1DNIddp.tex}
      \caption{Example of robust optimization using IAGO and the cost
        function $Q^{90\%}$. The function $f$ (dotted line), corrupted
        by a Gaussian noise on the factor (zero mean with a standard
        deviation of 0.2), is studied starting from an initial design
        of three points (as in Figure~\ref{fig:SUR1D}). Six additional
        evaluations are carried out (triangles), which are used to
        estimate the cost function based on the Kriging model (bold
        line), along with the estimated point mass distribution of the
        robust minimizers (\emph{right}). The cost function $Q^{90\%}$
        estimated, only for the sake of comparison, from the true
        function using Monte Carlo uncertainty propagation is also
        provided (mixed line).}
\label{fig:SUR1DROB}
  \end{figure}

Figure~\ref{fig:SUR1DROB} illustrates robust optimization using the
same function and initial design. The cost function used is the
90\%-quantile $Q^{90\%}$, which is computed on the surrogate model but
also, and only for the sake of comparison, on the true function using Monte
Carlo uncertainty propagation (the quantile is approximated using 5000
simulations). After six iterations of the robust optimization
algorithm, the distribution of the robust minimizer is sufficiently
peaked to give a good approximation of the true global robust
minimizer.

These result are encouraging as they show that the requirement of a
fast uncertainty reduction is met. The next section provides some more
examples, along with a comparison with EGO, the EI-based global
optimization algorithm.

\subsection{Empirical comparison with expected improvement}

\begin{figure}
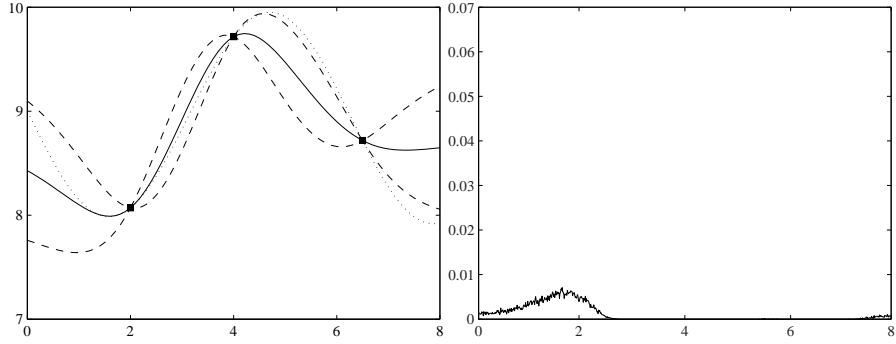
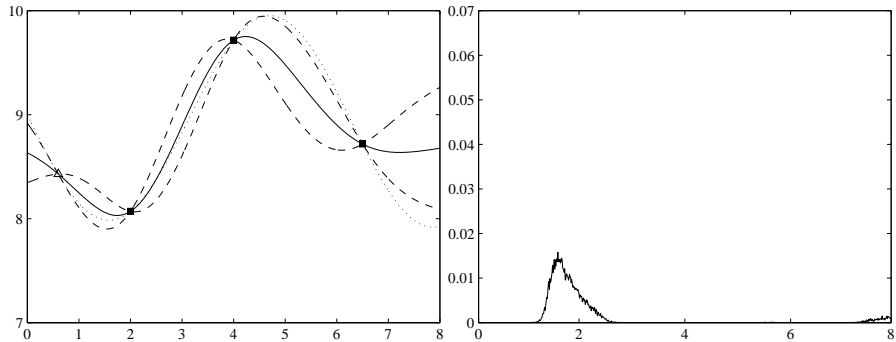
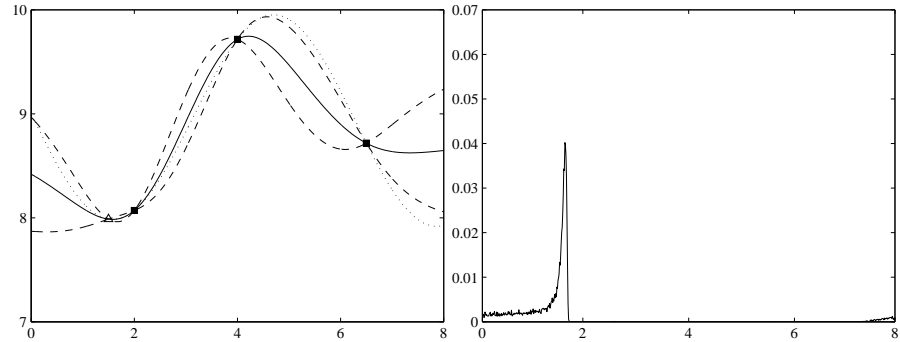

    \centering \subfigure[Initial prediction and density of the global
    minimizers]{\input{./courbes/SUR1initial.tex}
      \input{./courbes/SUR1initialddp.tex}} \subfigure[Prediction and
    density of the global minimizers after an additional evaluation of
    $f$ chosen with
    EI]{\label{fig:criterionEI}\input{./courbes/EI1.tex}
      \input{./courbes/EI1ddp.tex}} \subfigure[Prediction and density
    of the global minimizers after an additional evaluation of $f$
    chosen with conditional entropy]{\label{fig:criterionSUR}
      \input{./courbes/SUR1.tex} \input{./courbes/SUR1ddp.tex} }
    \caption{Comparison between conditional entropy and EI: the \emph{left}
      side contains the Kriging predictions before and after an
      additional evaluation chosen with either EI or conditional entropy, while
      the \emph{right} side presents the corresponding conditional
      density of the global minimizers.}
    \label{fig:criterion}
  \end{figure}

Consider first the function described by Figure~\ref{fig:criterion}.
Given an initial design of three points, both the EI and conditional
entropy are computed. Their optimization provides two candidate
evaluation point for $f$, which are also presented on
Figure~\ref{fig:criterion}, along with the post-evaluation prediction
and conditional point mass distribution for $\Minimizergrid$.  For
this example, the regularity parameter of the Mat\'ern covariance is set
a priori to a high value (2.5), and it is therefore likely that two
evaluations close to one another, as proposed by the IAGO algorithm,
will give some valuable information about the part of $\XX$ located on
the left of the evaluation point, while improving the characterization
of the global minimizer. By taking in account the covariance of $F$
through conditional simulations, the conditional entropy uses regularity to conclude faster. The resulting conditional
point mass distribution of the minimizers is then generally more peaked
using the IAGO algorithm than using the EGO algorithm (as illustrated
by Figure~\ref{fig:criterionSUR} and Figure~\ref{fig:criterionEI}).

\begin{figure}
    \centering
    \subfigure[15 iterations using EGO]{
      % This file is generated by the MATLAB m-file laprint.m. It can be included
% into LaTeX documents using the packages graphicx, color and psfrag.
% It is accompanied by a postscript file. A sample LaTeX file is:
%    \documentclass{article}\usepackage{graphicx,color,psfrag}
%    \begin{document}\input{EI15}\end{document}
% See http://www.mathworks.de/matlabcentral/fileexchange/loadFile.do?objectId=4638
% for recent versions of laprint.m.
%
% created by:           LaPrint version 3.16 (13.9.2004)
% created on:           01-Jun-2006 15:04:08
% eps bounding box:     15 cm x 10.7813 cm
% comment:              
%
\begin{psfrags}%
\psfragscanon%
%
% text strings:
\psfrag{s01}[t][t]{\color[rgb]{0,0,0}\setlength{\tabcolsep}{0pt}\begin{tabular}{c}\large$x_1$\end{tabular}}%
\psfrag{s02}[b][b]{\color[rgb]{0,0,0}\setlength{\tabcolsep}{0pt}\begin{tabular}{c}\large$x_2$\end{tabular}}%
%
% xticklabels:
\psfrag{x01}[t][t]{0}%
\psfrag{x02}[t][t]{0.1}%
\psfrag{x03}[t][t]{0.2}%
\psfrag{x04}[t][t]{0.3}%
\psfrag{x05}[t][t]{0.4}%
\psfrag{x06}[t][t]{0.5}%
\psfrag{x07}[t][t]{0.6}%
\psfrag{x08}[t][t]{0.7}%
\psfrag{x09}[t][t]{0.8}%
\psfrag{x10}[t][t]{0.9}%
\psfrag{x11}[t][t]{1}%
\psfrag{x12}[t][t]{0}%
\psfrag{x13}[t][t]{0.1}%
\psfrag{x14}[t][t]{0.2}%
\psfrag{x15}[t][t]{0.3}%
\psfrag{x16}[t][t]{0.4}%
\psfrag{x17}[t][t]{0.5}%
\psfrag{x18}[t][t]{0.6}%
\psfrag{x19}[t][t]{0.7}%
\psfrag{x20}[t][t]{0.8}%
\psfrag{x21}[t][t]{0.9}%
\psfrag{x22}[t][t]{1}%
\psfrag{x23}[t][t]{-5}%
\psfrag{x24}[t][t]{0}%
\psfrag{x25}[t][t]{5}%
\psfrag{x26}[t][t]{10}%
%
% yticklabels:
\psfrag{v01}[r][r]{0}%
\psfrag{v02}[r][r]{0.1}%
\psfrag{v03}[r][r]{0.2}%
\psfrag{v04}[r][r]{0.3}%
\psfrag{v05}[r][r]{0.4}%
\psfrag{v06}[r][r]{0.5}%
\psfrag{v07}[r][r]{0.6}%
\psfrag{v08}[r][r]{0.7}%
\psfrag{v09}[r][r]{0.8}%
\psfrag{v10}[r][r]{0.9}%
\psfrag{v11}[r][r]{1}%
\psfrag{v12}[r][r]{0}%
\psfrag{v13}[r][r]{0.1}%
\psfrag{v14}[r][r]{0.2}%
\psfrag{v15}[r][r]{0.3}%
\psfrag{v16}[r][r]{0.4}%
\psfrag{v17}[r][r]{0.5}%
\psfrag{v18}[r][r]{0.6}%
\psfrag{v19}[r][r]{0.7}%
\psfrag{v20}[r][r]{0.8}%
\psfrag{v21}[r][r]{0.9}%
\psfrag{v22}[r][r]{1}%
\psfrag{v23}[r][r]{0}%
\psfrag{v24}[r][r]{5}%
\psfrag{v25}[r][r]{10}%
\psfrag{v26}[r][r]{15}%
%
% Figure:
\resizebox{6cm}{!}{\includegraphics{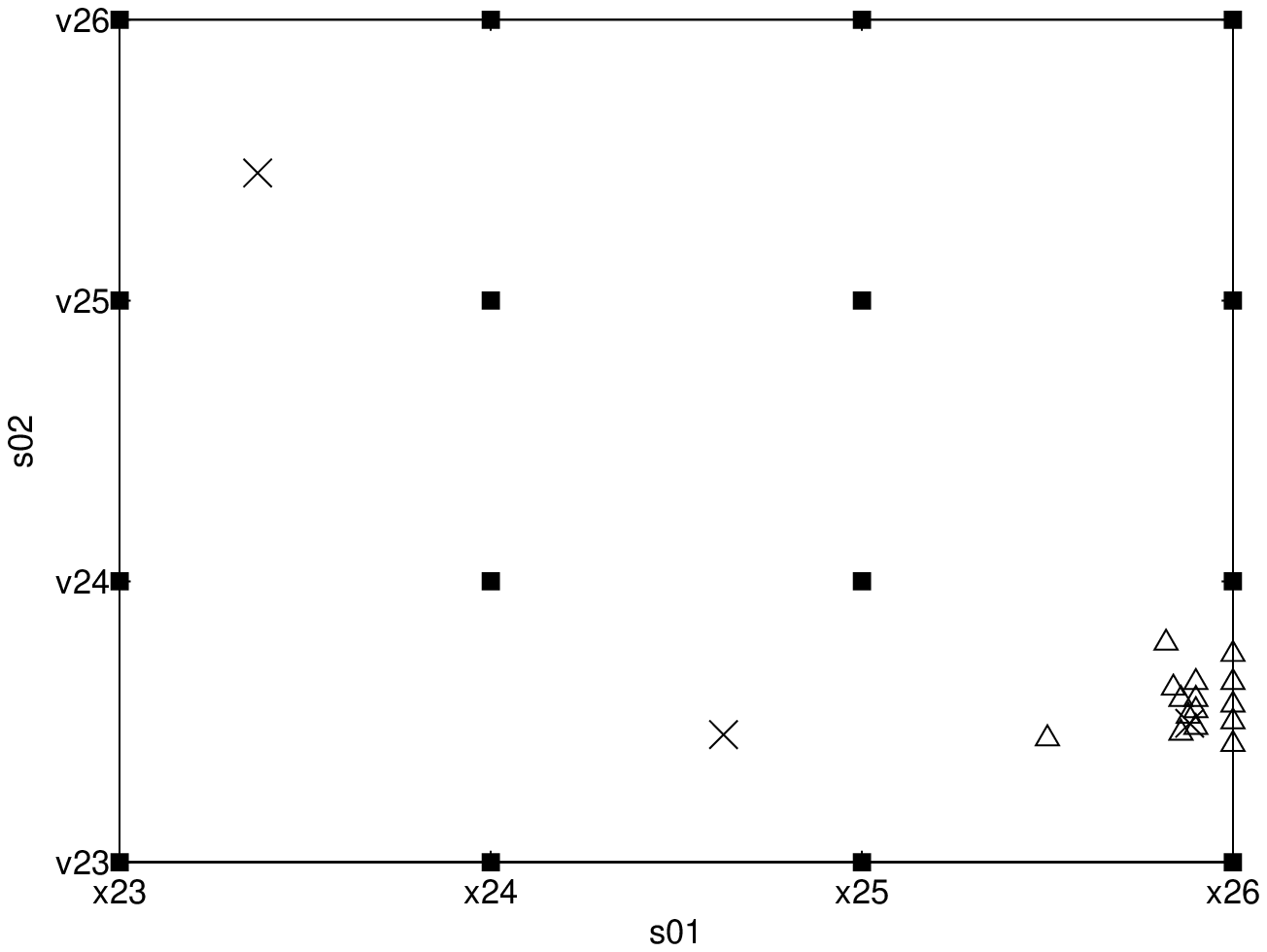}}%
\end{psfrags}%
%
% End EI15.tex
}
    \subfigure[35 iterations using EGO]{
      % This file is generated by the MATLAB m-file laprint.m. It can be included
% into LaTeX documents using the packages graphicx, color and psfrag.
% It is accompanied by a postscript file. A sample LaTeX file is:
%    \documentclass{article}\usepackage{graphicx,color,psfrag}
%    \begin{document}\input{EI2D}\end{document}
% See http://www.mathworks.de/matlabcentral/fileexchange/loadFile.do?objectId=4638
% for recent versions of laprint.m.
%
% created by:           LaPrint version 3.16 (13.9.2004)
% created on:           01-Jun-2006 15:02:51
% eps bounding box:     15 cm x 10.7813 cm
% comment:              
%
\begin{psfrags}%
\psfragscanon%
%
% text strings:
\psfrag{s01}[t][t]{\color[rgb]{0,0,0}\setlength{\tabcolsep}{0pt}\begin{tabular}{c}\large$x_1$\end{tabular}}%
\psfrag{s02}[b][b]{\color[rgb]{0,0,0}\setlength{\tabcolsep}{0pt}\begin{tabular}{c}\large$x_2$\end{tabular}}%
%
% xticklabels:
\psfrag{x01}[t][t]{0}%
\psfrag{x02}[t][t]{0.1}%
\psfrag{x03}[t][t]{0.2}%
\psfrag{x04}[t][t]{0.3}%
\psfrag{x05}[t][t]{0.4}%
\psfrag{x06}[t][t]{0.5}%
\psfrag{x07}[t][t]{0.6}%
\psfrag{x08}[t][t]{0.7}%
\psfrag{x09}[t][t]{0.8}%
\psfrag{x10}[t][t]{0.9}%
\psfrag{x11}[t][t]{1}%
\psfrag{x12}[t][t]{0}%
\psfrag{x13}[t][t]{0.1}%
\psfrag{x14}[t][t]{0.2}%
\psfrag{x15}[t][t]{0.3}%
\psfrag{x16}[t][t]{0.4}%
\psfrag{x17}[t][t]{0.5}%
\psfrag{x18}[t][t]{0.6}%
\psfrag{x19}[t][t]{0.7}%
\psfrag{x20}[t][t]{0.8}%
\psfrag{x21}[t][t]{0.9}%
\psfrag{x22}[t][t]{1}%
\psfrag{x23}[t][t]{-5}%
\psfrag{x24}[t][t]{0}%
\psfrag{x25}[t][t]{5}%
\psfrag{x26}[t][t]{10}%
%
% yticklabels:
\psfrag{v01}[r][r]{0}%
\psfrag{v02}[r][r]{0.1}%
\psfrag{v03}[r][r]{0.2}%
\psfrag{v04}[r][r]{0.3}%
\psfrag{v05}[r][r]{0.4}%
\psfrag{v06}[r][r]{0.5}%
\psfrag{v07}[r][r]{0.6}%
\psfrag{v08}[r][r]{0.7}%
\psfrag{v09}[r][r]{0.8}%
\psfrag{v10}[r][r]{0.9}%
\psfrag{v11}[r][r]{1}%
\psfrag{v12}[r][r]{0}%
\psfrag{v13}[r][r]{0.1}%
\psfrag{v14}[r][r]{0.2}%
\psfrag{v15}[r][r]{0.3}%
\psfrag{v16}[r][r]{0.4}%
\psfrag{v17}[r][r]{0.5}%
\psfrag{v18}[r][r]{0.6}%
\psfrag{v19}[r][r]{0.7}%
\psfrag{v20}[r][r]{0.8}%
\psfrag{v21}[r][r]{0.9}%
\psfrag{v22}[r][r]{1}%
\psfrag{v23}[r][r]{0}%
\psfrag{v24}[r][r]{5}%
\psfrag{v25}[r][r]{10}%
\psfrag{v26}[r][r]{15}%
%
% Figure:
\resizebox{6cm}{!}{\includegraphics{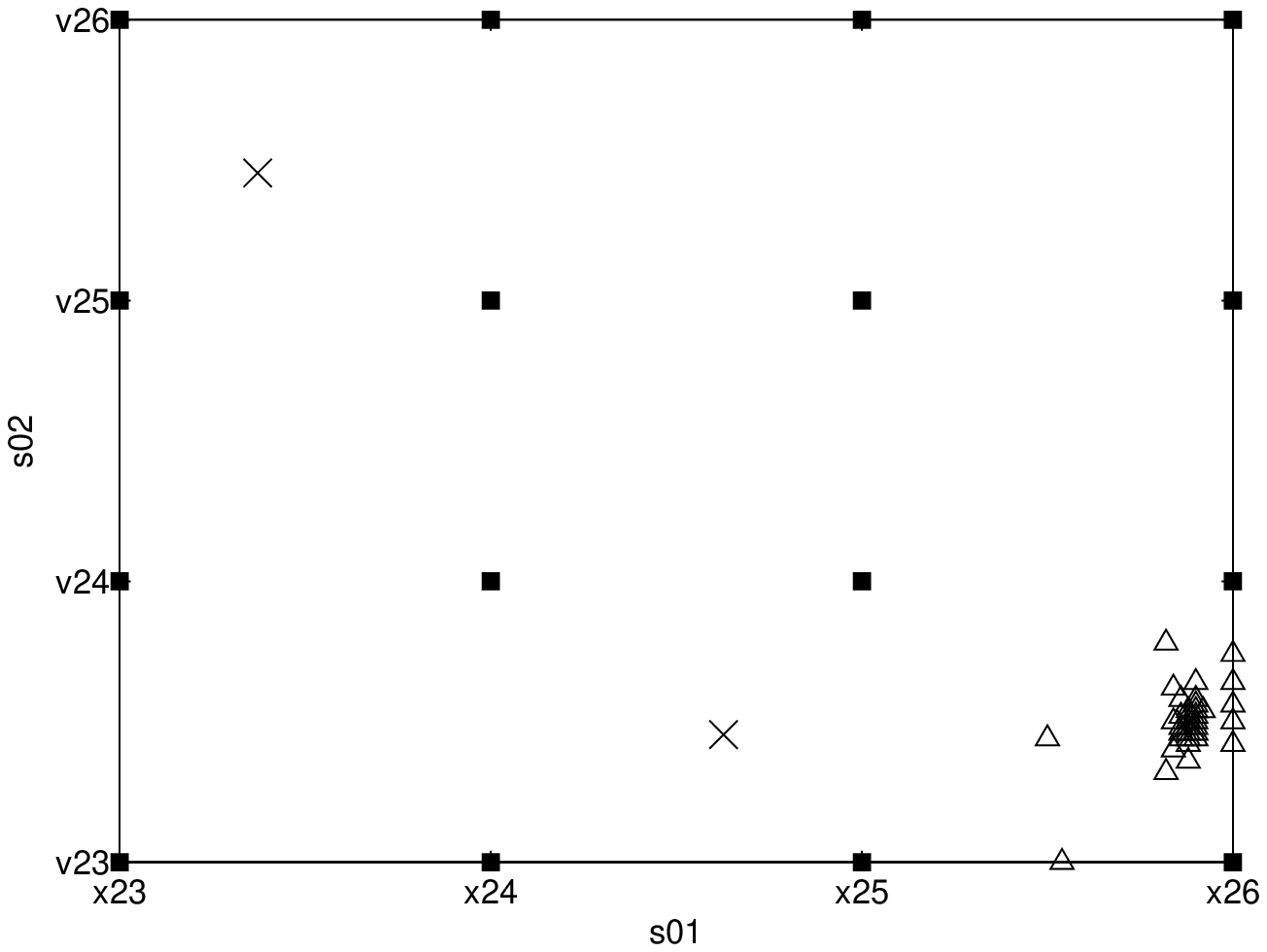}}%
\end{psfrags}%
%
% End EI2D.tex
}
    \subfigure[15 iterations using IAGO]{
      % This file is generated by the MATLAB m-file laprint.m. It can be included
% into LaTeX documents using the packages graphicx, color and psfrag.
% It is accompanied by a postscript file. A sample LaTeX file is:
%    \documentclass{article}\usepackage{graphicx,color,psfrag}
%    \begin{document}\input{SUR15}\end{document}
% See http://www.mathworks.de/matlabcentral/fileexchange/loadFile.do?objectId=4638
% for recent versions of laprint.m.
%
% created by:           LaPrint version 3.16 (13.9.2004)
% created on:           01-Jun-2006 15:05:24
% eps bounding box:     15 cm x 10.7813 cm
% comment:              
%
\begin{psfrags}%
\psfragscanon%
\large
%
% text strings:
\psfrag{s01}[l][l]{\color[rgb]{0,0,0}\setlength{\tabcolsep}{0pt}\begin{tabular}{l}\large 1\end{tabular}}%
\psfrag{s02}[l][l]{\color[rgb]{0,0,0}\setlength{\tabcolsep}{0pt}\begin{tabular}{l}\large 2\end{tabular}}%
\psfrag{s03}[l][l]{\color[rgb]{0,0,0}\setlength{\tabcolsep}{0pt}\begin{tabular}{l}\large 3\end{tabular}}%
\psfrag{s04}[l][l]{\color[rgb]{0,0,0}\setlength{\tabcolsep}{0pt}\begin{tabular}{l}\large 4\end{tabular}}%
\psfrag{s05}[l][l]{\color[rgb]{0,0,0}\setlength{\tabcolsep}{0pt}\begin{tabular}{l}\large 5\end{tabular}}%
\psfrag{s06}[l][l]{\color[rgb]{0,0,0}\setlength{\tabcolsep}{0pt}\begin{tabular}{l}\large 6\end{tabular}}%
\psfrag{s07}[l][l]{\color[rgb]{0,0,0}\setlength{\tabcolsep}{0pt}\begin{tabular}{l}\large 7\end{tabular}}%
\psfrag{s08}[l][l]{\color[rgb]{0,0,0}\setlength{\tabcolsep}{0pt}\begin{tabular}{l}\large 8\end{tabular}}%
\psfrag{s09}[l][l]{\color[rgb]{0,0,0}\setlength{\tabcolsep}{0pt}\begin{tabular}{l}\large 9\end{tabular}}%
\psfrag{s10}[l][l]{\color[rgb]{0,0,0}\setlength{\tabcolsep}{0pt}\begin{tabular}{l}\large10\end{tabular}}%
\psfrag{s11}[l][l]{\color[rgb]{0,0,0}\setlength{\tabcolsep}{0pt}\begin{tabular}{l}\large11\end{tabular}}%
\psfrag{s12}[l][l]{\color[rgb]{0,0,0}\setlength{\tabcolsep}{0pt}\begin{tabular}{l}\large12\end{tabular}}%
\psfrag{s13}[l][l]{\color[rgb]{0,0,0}\setlength{\tabcolsep}{0pt}\begin{tabular}{l}\large13\end{tabular}}%
\psfrag{s14}[l][l]{\color[rgb]{0,0,0}\setlength{\tabcolsep}{0pt}\begin{tabular}{l}\large14\end{tabular}}%
\psfrag{s15}[l][l]{\color[rgb]{0,0,0}\setlength{\tabcolsep}{0pt}\begin{tabular}{l}\large15\end{tabular}}%
\psfrag{s16}[l][l]{\color[rgb]{0,0,0}\setlength{\tabcolsep}{0pt}\begin{tabular}{l}\large16\end{tabular}}%
\psfrag{s17}[t][t]{\color[rgb]{0,0,0}\setlength{\tabcolsep}{0pt}\begin{tabular}{c}\large$x_1$\end{tabular}}%
\psfrag{s18}[b][b]{\color[rgb]{0,0,0}\setlength{\tabcolsep}{0pt}\begin{tabular}{c}\large$x_2$\end{tabular}}%
%
% xticklabels:
\psfrag{x01}[t][t]{0}%
\psfrag{x02}[t][t]{0.1}%
\psfrag{x03}[t][t]{0.2}%
\psfrag{x04}[t][t]{0.3}%
\psfrag{x05}[t][t]{0.4}%
\psfrag{x06}[t][t]{0.5}%
\psfrag{x07}[t][t]{0.6}%
\psfrag{x08}[t][t]{0.7}%
\psfrag{x09}[t][t]{0.8}%
\psfrag{x10}[t][t]{0.9}%
\psfrag{x11}[t][t]{1}%
\psfrag{x12}[t][t]{-5}%
\psfrag{x13}[t][t]{0}%
\psfrag{x14}[t][t]{5}%
\psfrag{x15}[t][t]{10}%
%
% yticklabels:
\psfrag{v01}[r][r]{0}%
\psfrag{v02}[r][r]{0.1}%
\psfrag{v03}[r][r]{0.2}%
\psfrag{v04}[r][r]{0.3}%
\psfrag{v05}[r][r]{0.4}%
\psfrag{v06}[r][r]{0.5}%
\psfrag{v07}[r][r]{0.6}%
\psfrag{v08}[r][r]{0.7}%
\psfrag{v09}[r][r]{0.8}%
\psfrag{v10}[r][r]{0.9}%
\psfrag{v11}[r][r]{1}%
\psfrag{v12}[r][r]{0}%
\psfrag{v13}[r][r]{5}%
\psfrag{v14}[r][r]{10}%
\psfrag{v15}[r][r]{15}%
%
% Figure:
\resizebox{6cm}{!}{\includegraphics{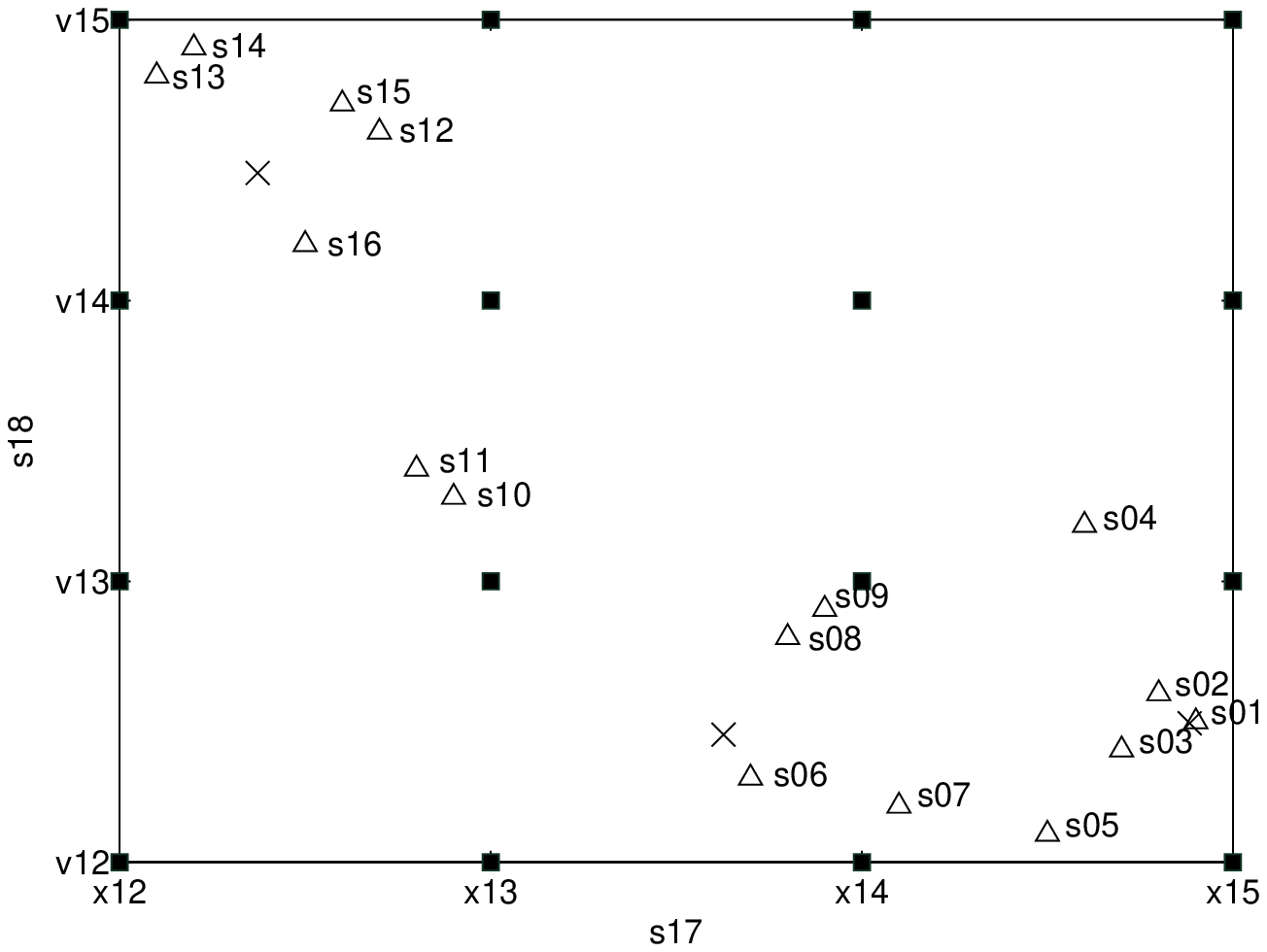}}%
\end{psfrags}%
%
% End SUR15.tex
}
    \subfigure[35 iterations using IAGO]{
      % This file is generated by the MATLAB m-file laprint.m. It can be included
% into LaTeX documents using the packages graphicx, color and psfrag.
% It is accompanied by a postscript file. A sample LaTeX file is:
%    \documentclass{article}\usepackage{graphicx,color,psfrag}
%    \begin{document}\input{SUR2D}\end{document}
% See http://www.mathworks.de/matlabcentral/fileexchange/loadFile.do?objectId=4638
% for recent versions of laprint.m.
%
% created by:           LaPrint version 3.16 (13.9.2004)
% created on:           01-Jun-2006 15:04:58
% eps bounding box:     15 cm x 10.7813 cm
% comment:              
%
\begin{psfrags}%
\psfragscanon%
%
% text strings:
\psfrag{s01}[t][t]{\color[rgb]{0,0,0}\setlength{\tabcolsep}{0pt}\begin{tabular}{c}\large$x_1$\end{tabular}}%
\psfrag{s02}[b][b]{\color[rgb]{0,0,0}\setlength{\tabcolsep}{0pt}\begin{tabular}{c}\large$x_2$\end{tabular}}%
%
% xticklabels:
\psfrag{x01}[t][t]{0}%
\psfrag{x02}[t][t]{0.1}%
\psfrag{x03}[t][t]{0.2}%
\psfrag{x04}[t][t]{0.3}%
\psfrag{x05}[t][t]{0.4}%
\psfrag{x06}[t][t]{0.5}%
\psfrag{x07}[t][t]{0.6}%
\psfrag{x08}[t][t]{0.7}%
\psfrag{x09}[t][t]{0.8}%
\psfrag{x10}[t][t]{0.9}%
\psfrag{x11}[t][t]{1}%
\psfrag{x12}[t][t]{0}%
\psfrag{x13}[t][t]{0.1}%
\psfrag{x14}[t][t]{0.2}%
\psfrag{x15}[t][t]{0.3}%
\psfrag{x16}[t][t]{0.4}%
\psfrag{x17}[t][t]{0.5}%
\psfrag{x18}[t][t]{0.6}%
\psfrag{x19}[t][t]{0.7}%
\psfrag{x20}[t][t]{0.8}%
\psfrag{x21}[t][t]{0.9}%
\psfrag{x22}[t][t]{1}%
\psfrag{x23}[t][t]{-5}%
\psfrag{x24}[t][t]{0}%
\psfrag{x25}[t][t]{5}%
\psfrag{x26}[t][t]{10}%
%
% yticklabels:
\psfrag{v01}[r][r]{0}%
\psfrag{v02}[r][r]{0.1}%
\psfrag{v03}[r][r]{0.2}%
\psfrag{v04}[r][r]{0.3}%
\psfrag{v05}[r][r]{0.4}%
\psfrag{v06}[r][r]{0.5}%
\psfrag{v07}[r][r]{0.6}%
\psfrag{v08}[r][r]{0.7}%
\psfrag{v09}[r][r]{0.8}%
\psfrag{v10}[r][r]{0.9}%
\psfrag{v11}[r][r]{1}%
\psfrag{v12}[r][r]{0}%
\psfrag{v13}[r][r]{0.1}%
\psfrag{v14}[r][r]{0.2}%
\psfrag{v15}[r][r]{0.3}%
\psfrag{v16}[r][r]{0.4}%
\psfrag{v17}[r][r]{0.5}%
\psfrag{v18}[r][r]{0.6}%
\psfrag{v19}[r][r]{0.7}%
\psfrag{v20}[r][r]{0.8}%
\psfrag{v21}[r][r]{0.9}%
\psfrag{v22}[r][r]{1}%
\psfrag{v23}[r][r]{0}%
\psfrag{v24}[r][r]{5}%
\psfrag{v25}[r][r]{10}%
\psfrag{v26}[r][r]{15}%
%
% Figure:
\resizebox{6cm}{!}{\includegraphics{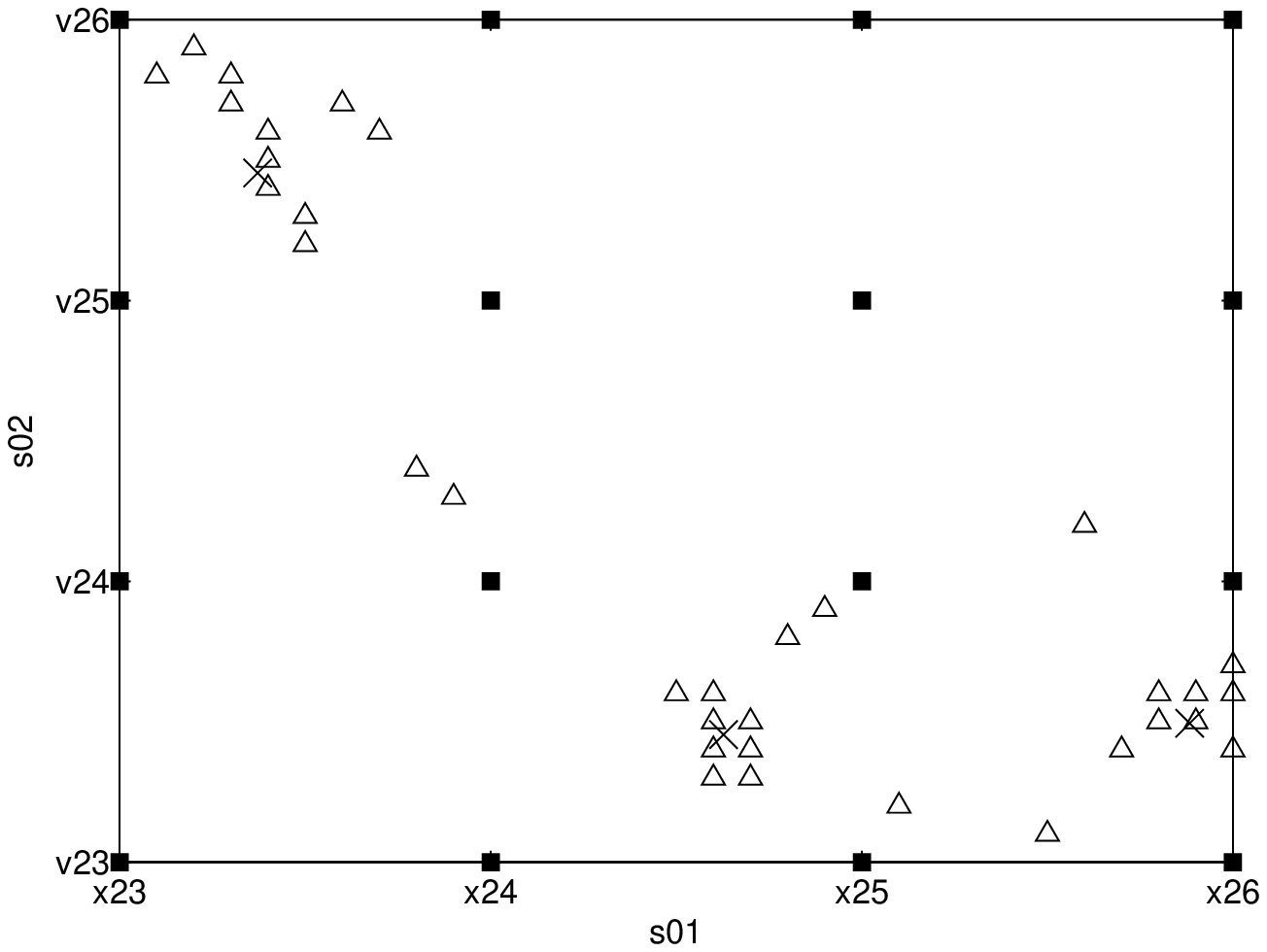}}%
\end{psfrags}%
%
% End SUR2D.tex
}
    \caption{fifteen iterations of two optimization algorithms, that
      differ by their criteria for selecting evaluation points for
      $f$, on the Branin function: (\emph{top}) the EI criterion  is
      used, (\emph{bottom}) the conditional entropy criterion is used with a
      thousand candidate evaluation points for $f$ set on a regular
      grid (squares account for initial data, triangles for new
      evaluations, and crosses give the actual locations of the three
      global minimizers).}
    \label{fig:Comparaison2D}
  \end{figure}

  Consider now the Branin function (see, for instance,
  (\cite{Dixon1978})), defined as
\begin{displaymath}
    \begin{array}{r r c l}
      f :   &  [-5,10] \times[0,15] &   \longrightarrow & \mathbb{R}  \\
      {}       &   (x_1,x_2) &   \longmapsto     & \left( x_2-\frac{5.1}{4\pi^2}x_1^2+\frac{5}{\pi}x_1-6\right)^2\\
      & & &   +10\left(1-\frac{1}{8\pi}\right)\cos(x_1)+10 \; .\\
    \end{array}
\end{displaymath}
It has three global minimizers $\minimizer_1\approx( -3.14,12.27 )\Tr$,
$\minimizer_2\approx (3.14,2.27)\Tr$ and $\minimizer_3\approx( 9.42,2.47 )\Tr$, and the
global minimum is approximately equal to 0.4. Given an initial uniform
design of fifteen points, fifteen additional points are iteratively
selected and evaluated using the IAGO and EGO algorithms. These
parameters are estimated on the initial design, and kept unchanged
during both procedures. The positions of these points are presented on
Figure~\ref{fig:Comparaison2D} (left), along with the three global
minimizers. Table \ref{tab:resultsBranin} summarizes the results
obtained with EGO and IAGO, based on the final Kriging models obtained
with both approaches. Note that the EI criterion in EGO is maximized
with a high precision, while the conditional entropy in IAGO is
computed over a thousand candidate evaluation points located on a
regular grid.
\begin{table}
\centering
\scriptsize
\begin{tabular}{|p{4 cm} |c |c |c |c|}
\hline
 & \multicolumn{2}{c|}{EGO} & \multicolumn{2}{c|}{IAGO} \\
& 15 iterations & 35 iterations & 15 iterations & 35 iterations \\
\hline
Euclidean distance between $ \minimizer_1$ and its final estimate& $3.22$&$3.22$ & $\bm{2.18}$&$\bm{0.23}$ \\
\hline
Value of the true function at estimated minimizer & $17.95$&$17.95$ & $\bm{2.59}$& $\bm{0.40}$\\
\hline
Euclidean distance between $ \minimizer_2$ and its final estimate & $2.40$& $2.40$& $\bm{0.44}$&$\bm{0.18}$\\
\hline
Value of the true function at estimated minimizer & $13.00$&$13.00$ & $\bm{0.85}$& $\bm{0.42}$ \\ 
\hline 
Euclidean distance between $ \minimizer_3$ and its final estimate & $\bm{0.04}$& $\bm{0.04}$& 0.82& 0.23\\
\hline
Value of the true function at estimated minimizer &$\bm{0.40}$& $\bm{0.40}$& 1.94&0.44 \\  
\hline
\end{tabular}
\caption[]{Estimation results for the Branin 
  function using evaluations of Figure~\ref{fig:Comparaison2D}}  
\label{tab:resultsBranin}
\end{table}
It appears nevertheless that the algorithm using EI stalls on a single
global minimizer, while the conditional entropy allows a relatively
fast estimation of all three of them. Besides, the search is more
global with IAGO, which yields a better global approximation of the
supposedly unknown function. If twenty additional evaluations are
carried out (as presented in the right part of
Figure~\ref{fig:Comparaison2D}), the final Kriging prediction using
the SUR approach estimates the minimum with an error of less than 0.05
for all three minimizers (cf.  Table~\ref{tab:resultsBranin}), while
the EI criterion does not improve the information on any minimizer any
further.  The difference between the two strategies is clearly
evidenced. The EI criterion, overestimating the confidence in the
initial prediction, has led to performing evaluations extremely close
to one another, for a very small information gain. In a context of
expensive function evaluation, this is highly detrimental. The entropy
criterion, using the same covariance parameters, does not stack points
almost at the same location before having identified the most likely
zones for the minimizers. The use of what has been assumed and learned
about the function is clearly more efficient in this case, and this
property should be highly attractive when dealing with problems of
higher dimension.

\section{Discussion}
\label{sec:discussion}

\subsection{Robustness to uncertainty on the covariance parameters}
\label{sec:robustn-estim-covar}

\cite{Jones2001} studies the potential of Kriging-based
  global optimization methods such as EGO. One of his most important
  conclusions, is that these methods ``\emph{can perform poorly if the
    initial sample is highly deceptive}''. An eloquent example is
  provided on page 373, where a sinus function is sampled using its
  own period, leading to a flat prediction over the domain, associated
  with a small prediction error.

This potential for deception is present throughout the IAGO
  procedure, and should not be ignored.  To overcome this difficulty,
  several methods have been proposed (see, e.g., the Enhanced Method 4
  in (\cite{Jones2001}) or (\cite{Gutmann2001})), which
  achieve some sort of robustness to an underestimation of the
  prediction error and more generally to a bad choice of covariance.
  They seem to perform better than classical algorithms, including
  EGO.

Comparing the IAGO approach to such methods is an interesting
  topic for future research. The issue considered here was to
  demonstrate the interest of the conditional entropy criterion, and
  we feel that this should be done independently from the rest of the
  procedure.

It is of course essential to make the IAGO algorithm robust to
  errors in the estimation of the covariance parameters. In many
  industrial problems, this can be easily done by using prior
  knowledge on the unknown function to restrict the possible values
  for these parameters. For example, experts of the field often have
  information regarding the range of values attainable by the unknown
  function.  This information can be directly used to restrict the
  search space for the variance of the modeling process $F$, or even
  to choose it beforehand.

More generally, given the probabilistic framework used here, it
  should be relatively easy to develop a Bayesian or minimax extension
  of the IAGO algorithm to guide the estimation of the parameters of
  the covariance. A comparison with robust methods such as those detailed in
  (\cite{Jones2001}) will then be essential.

\subsection{Conclusions and perspectives}
\label{sec:concl-persp}

In this paper, a stepwise uncertainty reduction strategy has been used
for the sequential global optimization of expensive-to-evaluate
functions. This strategy iteratively selects a minimizer of the
conditional minimizer entropy as the new evaluation point $\x_{\mathrm{new}}$. To compute
this entropy, a Gaussian random model of the function evaluations is
used and the minimizer entropy is estimated through Kriging and
conditional simulations. At each iteration, the result of the
evaluation  at $\x_{\mathrm{new}}$ is incorporated in the data
base used to re-build the Kriging model (with a possible re-estimation
of the parameters of its covariance).

We have shown on some simple examples that, compared to the classical
EI-based algorithm EGO, the method proposed significantly reduces the
evaluation effort in the search for global optimizers. The SUR
strategy allows the optimization method to adapt the type of search to
the information available on the function. In particular, the
conditional entropy criterion makes full use of the assumed regularity
of the unknown function to balance global and local searches.

Choosing an adequate set of candidate points is a crucial point, as it
must allow a good estimation of a global minimizer of the criterion,
while keeping computation feasible. Promising results have already
been obtained with space-filling designs, and adaptive sampling based
on the conditional density of the global minimizers should provide
useful results as dimension increases.

Extension to constrained optimization is an obviously important
  topic of future investigations. When it is easy to discard the
  candidate points in $\XX$ that do not satisfy the constraints, the
  extension is trivial. For expensive-to-evaluate constraints, the
  extension is a major challenge.

Finally, the SUR strategy associated with conditioning by Kriging is a
promising solution for the robust optimization of expensive-to-evaluate
functions, a problem that is central to many industrial situations,
for which an efficient product design must be found in the presence of
significant uncertainty on the values actually taken by some factors
in mass production. In addition, robustness to the uncertainty
associated with the estimation of the parameters of the covariance
should also be sought.

\section{Appendix: modeling with Gaussian processes}
\label{sec:annex-model-with}
This section recalls the main concepts used in this paper, namely
Gaussian process modeling and Kriging. The major results will be
presented along with the general framework for the estimation of the
model parameters.

\subsection{Kriging when $f$ is evaluated exactly}
\label{sec:generalities}
Kriging (\cite{Matheron1963,Chiles1999}) is a prediction method based
on random processes that can be used to approximate or interpolate
data. It can also be understood as a kernel regression method, such as
\emph{splines} (\cite{Wahba1998}) or \emph{Support Vector Regression}
(\cite{Smola1998}). It originates from geostatistics and is
widely used in this domain since the 60s. Kriging is also known as the
\emph{Best Linear Unbiased Prediction} (BLUP) in statistics, and has
been more recently designated as Gaussian Processes (GP) in the 90s in
the machine learning community.

As mentioned in Section~\ref{sec:gauss-proc-model}, it is assumed that
the function $f$ is a sample path of a second-order Gaussian random
process $F$. Denote by $m(\x)=E[F(\x)]$ the mean of $F(\x)$ and by
$k(\x,\y)$ its covariance, written as
$$k(\x,\y)=\mathrm{cov}(F(\x),F(\bm{y})).$$ 
Kriging then computes the BLUP of $F(\x)$, denoted by $\hat F(\x)$, in
the vector space generated by the evaluations
$\HH_\Sample=\mathrm{span}\{F(\x_1),\ldots,F(\x_n)\}$. As an element of
$\HH_\Sample$, $\pred(\x)$ can be written as
\begin{equation}
\label{eq:prediction2}
   \pred(\x) = \lambdas(\x)\Tr\bm{F}_\Sample\,,
\end{equation}
As the BLUP, $\pred(\x)$ must have the smallest variance for the
prediction error
\begin{equation}
\predvar^2(\x)=\expect[(\pred(\x)-F(\x))^2],
\label{eq:MSEsimple2}
\end{equation}
among all unbiased predictors. The variance of the prediction error
satisfies
\begin{equation}
  \label{eq:MSE}
    \predvar^2(\x)=k(\x,\x)+\lambdas(\x)\Tr \covmat \lambdas(\x)
  - 2\lambdas(\x)\Tr \bm{k}(\x),
\end{equation}
with 
$$
\covmat=\left(k(\x_i,\x_j)\right),\; (i,j) \in \llbracket 1,n \rrbracket^2
$$
the $\nobs\times\nobs$ covariance matrix of $F$ at evaluation points in
$\Sample$, and
$$
\bm{k}(\x)=[k(\x_1,\x),\ldots,k(\x_{\nobs},\x)]\Tr
$$
the vector of covariances between $F(\x)$ and $\bm{F}_\Sample$

The prediction method (\cite{Matheron1969}) assumes that the mean of
$F(\x)$ can
be written as a finite linear combination 
$$
  m(\x)=\bm{\beta}\Tr\bm{p}(\x),
$$
where $\bm{\beta}$ is a vector of fixed but unknown coefficients, and
$$\bm{p}(\x)=[p_1(\x), \ldots, p_l(\x)]\Tr$$ is a vector of known functions
of the factor vector $\x$.  Usually these functions are monomials of low
degree in the components of $\x$ (in practice, the degree does not
exceed two). These functions may be used to reflect some prior knowledge on
the unknown function. As we have none for the examples considered here,
we simply use an unknown constant.

The Kriging predictor at $\x$ is then the best linear predictor
subject to the unbiasedness constraint $\expect(\pred(\x))=m(\x)$,
whatever the unknown $\bm{\beta}$. The unbiasedness constraint translates
into
\begin{equation}
 \bm{\beta}\Tr\bm{P}\Tr\lambdas(\x)=\bm{\beta}\Tr\bm{p}(\x),
\label{eq:unbiaisdness}
\end{equation}
with
$$ \bm{ P}=\left( \begin{array}{c} \bm{p}(\x_1)\Tr\\
     \vdots\\\bm{p}(\x_\nobs)\Tr \end{array}\right).
$$
For (\ref{eq:unbiaisdness}) to be satisfied for all $\bm{\beta}$, the
Kriging coefficients must satisfy the linear constraints
\begin{equation}
  \label{eq:unbiased2}
 \bm{P}\Tr\lambdas(\x)=\bm{p}(\x),
\end{equation}
called \emph{universality constraints} by Matheron. At this point,
Kriging can be reformulated as follows: find the vector of Kriging
coefficients that minimizes the variance of the prediction error
(\ref{eq:MSE}) subject to the constraints (\ref{eq:unbiased2}). This
problem can be solved via a Lagrangian formulation, with $\bm{\mu}(\x)$
a vector of $l$~Lagrange multipliers.  The coefficients
$\lambdas(\x)$ are then solutions of the linear system of equations
\begin{equation}
\left( \begin{array}{c c} \covmat & \bm{P} \\ \bm{P}\Tr & \bm{0} \end{array} \right)
\left( \begin{array}{c} \lambdas(\x) \\  \bm{\mu}(\x) \end{array} \right)
= \left( \begin{array}{cc} \bm{k}(\x) \\ \bm{p}(\x) \end{array} \right) ,
\label{eq:Kr1}
\end{equation}
with $\bm{0}$ a matrix of zeros. A convenient expression for the
variance of the prediction error is obtained by substituting
$\bm{k}(\x)-\bm{P}\bm{\mu}(\x)$ for $\covmat\lambdas(\x)$ in
(\ref{eq:MSE}) as justified by (\ref{eq:Kr1}), to get
\begin{equation}
  \predvar^2(\x)=\expect
  \left[F(\x)-\pred(\x)\right]^2=k(\x,\x)-\lambdas(\x)\Tr\bm{k}(\x)-\bm{p}(\x)\Tr \bm{\mu}(\x)\: .
\label{eq:MSEKr}
\end{equation}
The variance of the prediction error at $\x$ can thus be computed
without any evaluation of $f$, using (\ref{eq:Kr1}) and
(\ref{eq:MSEKr}). It provides a measure of the quality associated with
the Kriging prediction. Evaluations of $f$ remain needed to
estimate the parameters of the covariance of $F$ (if any), as will be seen in
Section~\ref{sec:estim-param}.

Once $f$ has been evaluated at all evaluation points, the prediction
of the value taken by $f$ at $\x$ becomes 
\begin{equation}
  \predmean (\x)=\lambdas(\x)\Tr\bm{f}_\Sample\: ,
\label{eq:predmean}
\end{equation}
with $\bm{f}_{\Sample}=[f(\x_1),\ldots,f(\x_n)]\Tr$ ($\bm{f}_{\Sample}$ is
viewed as a sample value of $\bm{F}_{\Sample}$).

It is easy to check that (\ref{eq:Kr1}) implies that
$$
    \forall \:\x_i \: \in \: \Sample, \: \: \pred(\x_i)=F(\x_i). 
$$
The prediction of $f$ at $\x_i \:\in\: \Sample$ is then $f(\x_i)$, so Kriging is an
interpolation with the considerable advantage that it also accounts
for model uncertainty through an explicit characterization of the
prediction error.

{\bf Remark}: The Bayesian framework (see, for instance,
\cite{Williams1996}) is an alternative approach to derive the BLUP, in
which the model $F$ is viewed as a Bayesian prior on the output. In
the case of a zero-mean model, the conditional distribution of the
function is then Gaussian with mean
\begin{equation}
    \expect \left[F(\x) |\,\bm{F}_\Sample=\bm{f}_\Sample \right]=
  \bm{k}(\x)\Tr\covmat^{-1} \bm{f}_\Sample,
  \label{eq:conditional mean}
\end{equation}
and variance
$$
 \mathrm{Var} \left[F(\x) |\,\bm{F}_\Sample=\bm{f}_\Sample \right] =
k(\x,\x)-\bm{k}(\x)\Tr\covmat^{-1} \bm{k}(\x),
$$
which are exactly the mean~(\ref{eq:predmean}) and
variance~(\ref{eq:MSEKr}) of the Kriging predictor for a model $F$
with zero mean. The Kriging predictor can also be viewed as the
conditional mean of $F(\x)$ in the case of an unknown mean, if the
universality constraints are viewed as a non-informative prior on
$\bm{\beta}$.

\subsection{Kriging when $f$ is evaluated approximately} 
\label{sec:dealing-with-noisy}
The Kriging predictor was previously defined as the element of the
space $\HH_\Sample$ generated by the random variables $F(\x_i)$ that
minimizes the prediction error. A natural step is to extend this
formulation to the case of a function whose evaluations are corrupted
by additive independent and identically distributed Gaussian noise
variables $\varepsilon_i$ with zero mean and variance $\sigma_{\varepsilon}^2$. The model of
the observations then becomes $F^{\rm obs}_{\x_i}=F(\x_i)+\varepsilon_i \:
i=1,\ldots, \nobs$, and the Kriging predictor for $F(\x)$ takes the form
$\pred(\x)=\lambdas(\x)\Tr\bm{F}_{\Sample}^{\rm obs}$ with
$\bm{F}_\Sample^{\rm obs}=\left[F^{\rm obs}_{\x1},\ldots,F^{\rm
    obs}_{\x_{\nobs}} \right]\Tr$. The unbiasedness constraint
(\ref{eq:unbiased2}) remain unchanged, while the mean-square error
(\ref{eq:MSEsimple}) becomes
$$
\expect[\pred(\bm{x})-F(\bm{x})]^2= k(\x,\x)+\lambdas(\x)\Tr (\covmat+
\sigma_{\varepsilon}^2\bm{I}_{\nobs}) \lambdas(\x)
- 2\lambdas(\x)\Tr \bm{k}(\x),
$$
with $\bm{I}_n$ the identity matrix. Finally, using Lagrange
multipliers as before, it is easy to show that the coefficients
$\lambdas(\x)$ of the prediction must satisfy
\begin{equation}
\left( \begin{array}{c c} \covmat+ \sigma_{\varepsilon}^2\bm{I}_{\nobs}& \bm{P} \\ \bm{P}\Tr & 0 \end{array} \right)
\left( \begin{array}{c} \lambdas(\x) \\  \bm{\mu}(\x) \end{array} \right)
= \left( \begin{array}{cc} \bm{k}(\x) \\ \bm{p}(\x) \end{array} \right) .
\end{equation}
The resulting prediction is no longer interpolative, but can still be
viewed as the mean of the conditional distribution of $F$. The
variance of the prediction error is again obtained using
(\ref{eq:MSEKr}).
\subsection{Covariance choice}
\label{sec:implementation}
Choosing a suitable covariance function $k(\cdot ,\cdot)$ for a given $f$ is a
recurrent and fundamental question. It involves the choice of a
parametrized class (or model) of covariance, and the estimation of its
parameters.
\subsubsection{Covariance classes}
The asymptotic theory of Kriging (\cite{Stein1999}) stresses the
importance of the behaviour of the covariance near the origin.  This
behaviour is indeed linked with the quadratic-mean regularity of the
random process.  For instance, if the covariance is continuous at the
origin, then the process will be continuous in quadratic mean. In
practice, one often uses covariances that are \emph{invariant by
  translation} (or equivalently \emph{stationary}), \emph{isotropic},
and such that regularity can be adjusted.  Non-stationary covariances
are seldom used in practice, as they make parameter estimation
particularly difficult (\cite{Chiles1999}).  Isotropy, however, is
not required and can even be inappropriate when the factors are of
different natures. An example of an anisotropic, stationary covariance
class is $k(\x,\y)=k(h)$, with
$h=\sqrt{(\bm{x}-\bm{y})\Tr\bm{A}(\bm{x}-\bm{y})}$ where
$(\bm{x},\bm{y})\in \XX^2$ and $\bm{A}$ is a symmetric positive definite
matrix.

\begin{figure}[tbp]
  \centering
  % This file is generated by the MATLAB m-file laprint.m. It can be included
% into LaTeX documents using the packages graphicx, color and psfrag.
% It is accompanied by a postscript file. A sample LaTeX file is:
%    \documentclass{article}\usepackage{graphicx,color,psfrag}
%    \begin{document}\input{covreg}\end{document}
% See http://www.mathworks.de/matlabcentral/fileexchange/loadFile.do?objectId=4638
% for recent versions of laprint.m.
%
% created by:           LaPrint version 3.16 (13.9.2004)
% created on:           01-Jun-2006 13:35:03
% eps bounding box:     15 cm x 10.7813 cm
% comment:              
%
\begin{psfrags}%
\psfragscanon%
%
% text strings:
\psfrag{s05}[t][t]{\color[rgb]{0,0,0}\setlength{\tabcolsep}{0pt}\begin{tabular}{c}$h$\end{tabular}}%
\psfrag{s06}[b][b]{\color[rgb]{0,0,0}\setlength{\tabcolsep}{0pt}\begin{tabular}{c}$k(h)$\end{tabular}}%
%
% xticklabels:
\psfrag{x01}[t][t]{0}%
\psfrag{x02}[t][t]{0.1}%
\psfrag{x03}[t][t]{0.2}%
\psfrag{x04}[t][t]{0.3}%
\psfrag{x05}[t][t]{0.4}%
\psfrag{x06}[t][t]{0.5}%
\psfrag{x07}[t][t]{0.6}%
\psfrag{x08}[t][t]{0.7}%
\psfrag{x09}[t][t]{0.8}%
\psfrag{x10}[t][t]{0.9}%
\psfrag{x11}[t][t]{1}%
\psfrag{x12}[t][t]{-2}%
\psfrag{x13}[t][t]{-1.5}%
\psfrag{x14}[t][t]{-1}%
\psfrag{x15}[t][t]{-0.5}%
\psfrag{x16}[t][t]{0}%
\psfrag{x17}[t][t]{0.5}%
\psfrag{x18}[t][t]{1}%
\psfrag{x19}[t][t]{1.5}%
\psfrag{x20}[t][t]{2}%
%
% yticklabels:
\psfrag{v01}[r][r]{0}%
\psfrag{v02}[r][r]{0.1}%
\psfrag{v03}[r][r]{0.2}%
\psfrag{v04}[r][r]{0.3}%
\psfrag{v05}[r][r]{0.4}%
\psfrag{v06}[r][r]{0.5}%
\psfrag{v07}[r][r]{0.6}%
\psfrag{v08}[r][r]{0.7}%
\psfrag{v09}[r][r]{0.8}%
\psfrag{v10}[r][r]{0.9}%
\psfrag{v11}[r][r]{1}%
\psfrag{v12}[r][r]{0}%
\psfrag{v13}[r][r]{0.1}%
\psfrag{v14}[r][r]{0.2}%
\psfrag{v15}[r][r]{0.3}%
\psfrag{v16}[r][r]{0.4}%
\psfrag{v17}[r][r]{0.5}%
\psfrag{v18}[r][r]{0.6}%
\psfrag{v19}[r][r]{0.7}%
\psfrag{v20}[r][r]{0.8}%
\psfrag{v21}[r][r]{0.9}%
\psfrag{v22}[r][r]{1}%
%
% Figure:
\resizebox{12cm}{!}{\includegraphics{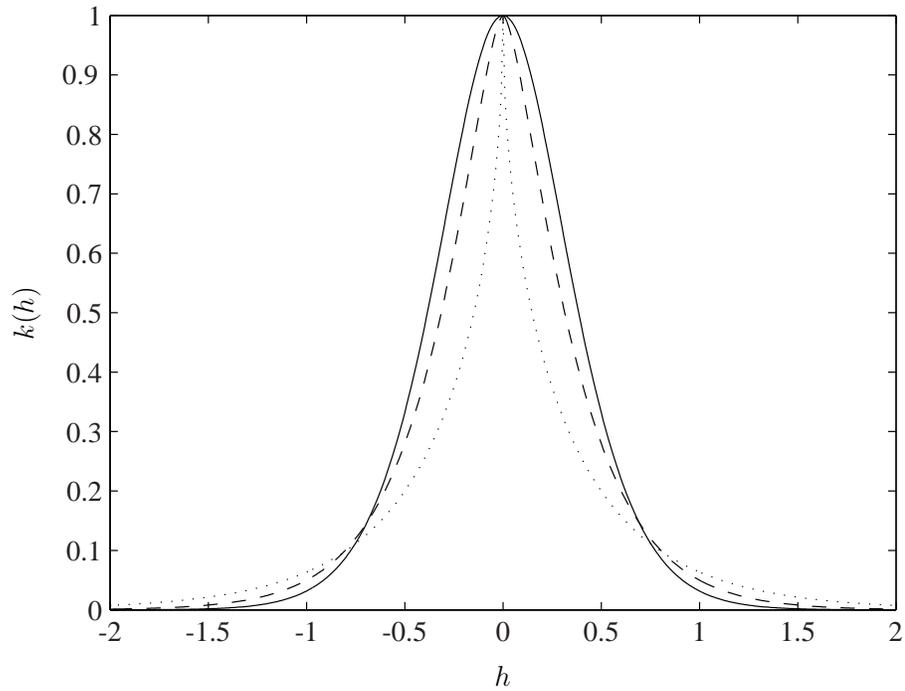}}%
\end{psfrags}%
%
% End covreg.tex

  \caption{Mat\'ern covariances for the
    parameterization  (\ref{eq:Matern})  with  $\rho=0.5$,  $\sigma^2=1$.  
    Solid  line corresponds  to $\nu=4$,  dashed line  to  $\nu=1$ and
    dotted line to $\nu=0.25$.}
  \label{fig:matern}
\end{figure}

\begin{figure}
  \centering
  \input{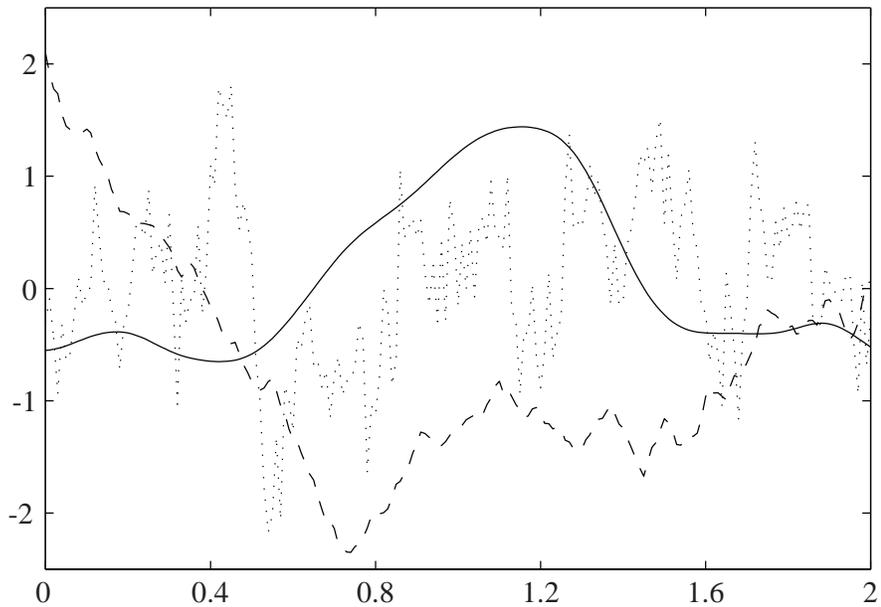}
  \caption{Three sample paths of a zero-mean Gaussian process with a
    Mat\'ern covariance. Conventions are as in Figure~\ref{fig:matern}:
    $\nu=4$ for the solid line, $\nu=1$ for the dashed line and $\nu=0.25$
    for the dotted line.}
  \label{fig:TrajMatern}
\end{figure}

A number of covariance classes are classically used (for instance
exponential $h \mapsto \sigma^2 \exp(-\theta |h|^\alpha )$, product of exponentials, or
polynomial). The \emph{Mat\'ern covariance} class offers the possibility
to adjust regularity with a single parameter (\cite{Stein1999}).
The Fourier transform of a Mat\'ern covariance is
\begin{equation}
    \hat{k}(\omega) = \frac{\sigma^2}{(\omega_0^2 +
  \omega^2)^{\nu+1/2}}\,, \quad \omega \in \mathbb{R}\,,
\label{eq:Matern}
\end{equation}
where $\nu$ controls the decay of $\hat{k}(\omega)$ at infinity and therefore
the regularity of the covariance at the origin. \cite{Stein1999})
advocates the use of the following parametrization of the Mat\'ern
class:
\begin{equation}
    k(h) =  \frac{\sigma^2}{2^{\nu-1}\Gamma(\nu)}\left(\frac{2\nu^{1/2}h}{\rho}\right)^\nu
  \mathcal{K}_\nu\left(\frac{2\nu^{1/2}h}{\rho}\right)\,,
\end{equation}
where $\mathcal{K}_\nu$ is the modified Bessel function of the second
kind (\cite{Yaglom1986}).  This parameterization is easy to
interpret, as $\nu$ controls regularity, $\sigma^2$ is the variance ($k(0) =
\sigma^2$), and $\rho$ represents the \emph{range} of the covariance,
\emph{i.e.}, the characteristic correlation distance. To stress the
significance and relevance of the regularity parameter,
Figure~\ref{fig:matern} shows the influence of $\nu$ on the covariance,
and Figure~\ref{fig:TrajMatern} demonstrates its impact on the sample
paths. Since Kriging assumes that $f$ is a sample path of $F$, a
careful choice of the parameters of the covariance is essential.
\subsubsection{Covariance parameters}
\label{sec:estim-param}
The parameters for a given covariance class can either be fixed using
prior knowledge on the system, or be estimated from experimental data.
In geostatistics, estimation is carried out using the adequacy between
the empirical and model covariances (\cite{Chiles1999}). In other
areas, cross validation (\cite{Wahba1998}) and maximum likelihood
(\cite{Stein1999}) are mostly employed. For simplicity and generality
reasons (\cite{Stein1999}), the maximum-likelihood method is preferred
here.

The maximum-likelihood estimate of the parameters of the covariance
maximizes the probability density of the data.  Using the joint
probability density of the observed Gaussian vector, and assuming that
the mean of $F(\x)$ is zero for the sake of simplicity, the
maximum-likelihood estimate of the vector $\bm{\theta}$ of the
covariance parameters is obtained (see, for instance, \cite{Vecchia1988}) by
minimizing the negative log-likelihood
\begin{equation}
   l(\param) = \frac{\nobs}{2} \log 2\pi + \frac{1}{2} \log \det {\covmat(\param)} +
  \frac{1}{2} \bm{f}_\Sample\Tr \covmat(\param)^{-1} \bm{f}_\Sample\,, 
\end{equation}
In the case of an unknown mean for $F(\x)$, it is possible to estimate
the parameters, using for example the \emph{REstricted Maximum
 Likelihood} (REML, see \cite{Stein1999}). This is the
 approach used for the examples in this paper.

\begin{figure}
  \centering \input{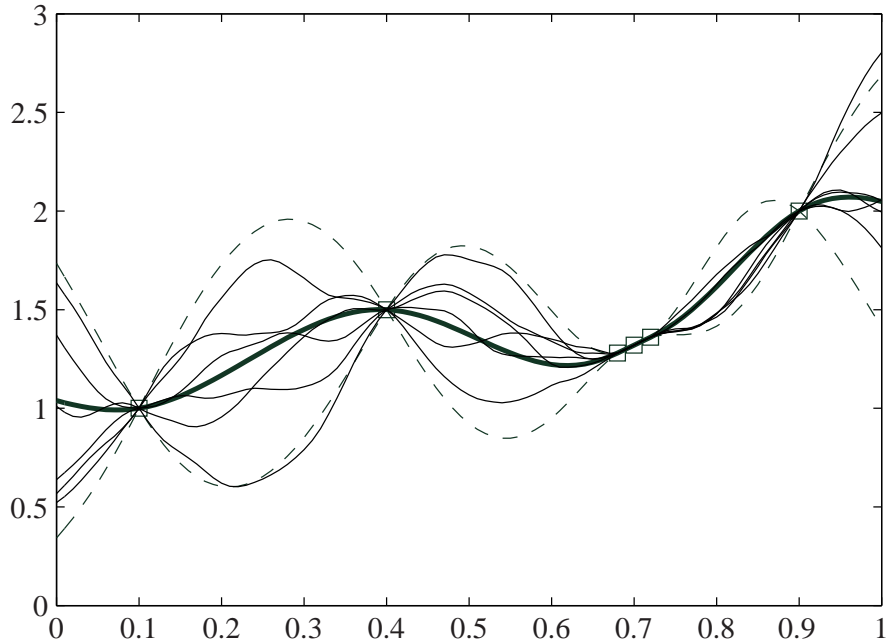}
  \caption{Example of Kriging interpolation (bold line) for a function
    of one variable. The data are represented by squares, and the
    covariance parameters were estimated by REML. Dashed
    lines delimit 95\% confidence region for the prediction. The
    thin solid lines are examples of conditional simulations.}
  \label{fig:ExempleKr}
\end{figure}

Figure~\ref{fig:ExempleKr} illustrates prediction by Kriging with a
Mat\'ern covariance, the parameters of which have been estimated by
REML. The prediction interpolates the data, and confidence intervals
are deduced from the square root of the variance of the prediction
error to assess the quality of the prediction between data.
Figure~\ref{fig:ExempleKr} also contains a series of conditional
simulations (obtained with the method explained in
Section~\ref{sec:conditionning-by-Kriging}), namely sample paths of
$F$ that interpolate the data. As implied by (\ref{eq:conditional
  mean}), the Kriging prediction is the mean of these conditional simulations.

\section*{Acknowledgements}
The authors wish to thank Donald R. Jones for his comments that
greatly contributed to improving the accuracy and clarity of this
paper.

%
%
%
%
%%%%%%%%%%%%%%%%%%%%%% 
%                    %      
%                    % 
%    Bibliography    % 
%                    % 
%                    %
%%%%%%%%%%%%%%%%%%%%%%

 \bibliographystyle{plainnat}
 \bibliography{./sur}

\end{document}